\documentclass[11pt,a4paper]{article}
\usepackage{amsfonts}
\usepackage{latexsym}
\usepackage{amsmath}
\usepackage{amssymb}
\usepackage{amsthm}
\usepackage{amsmath}
\usepackage{verbatim}
\usepackage{hyperref}
\usepackage{enumerate}

\usepackage[T1]{fontenc}
\usepackage{amsfonts}
\usepackage{graphicx}
\usepackage{epstopdf}
\ifpdf
  \DeclareGraphicsExtensions{.eps,.pdf,.png,.jpg}
\else
  \DeclareGraphicsExtensions{.eps}
\fi

\usepackage{tikz}
\usetikzlibrary{calc} 
\usepackage{array}  
\usepackage{amssymb} 
\usepackage[shortlabels]{enumitem} 

\usepackage{amsopn}

\newtheorem{theorem}{Theorem}[section]

\newtheorem{lemma}[theorem]{Lemma}
\newtheorem{proposition}[theorem]{Proposition}
\newtheorem{conjecture}[theorem]{Conjecture}
\newtheorem{definition}[theorem]{Definition}
\theoremstyle{definition}
\newtheorem{remark}[theorem]{\textbf{Remark}}
\theoremstyle{proof}

\numberwithin{equation}{section}

\numberwithin{figure}{section}

\newcommand{\mrm}{\mathrm}

\newcommand{\vol}{\mrm{vol}}
\newcommand{\Vol}{V_\mu}
\newcommand{\Volm}{V}
\newcommand{\Volg}{V_{\gamma^n}}  
\newcommand{\Area}{A_\mu}
\newcommand{\Ric}{\mrm{Ric}}
\newcommand{\tr}{\mrm{tr}}

\DeclareMathOperator*{\argmin}{arg\,min} 

\newcommand{\n}{\mathbf{n}}
\renewcommand{\n}{\mathfrak{n}}
\renewcommand{\c}{\mathbf{c}}
\renewcommand{\k}{\mathbf{k}}

\newcommand{\G}{\mathbb{G}}
\newcommand{\HH}{\mathbb{H}}
\renewcommand{\SS}{\mathcal{S}}
\newcommand{\Y}{\mathbf{Y}}
\newcommand{\T}{\mathbf{T}}
\newcommand{\F}{\mathcal{F}}
\newcommand{\FF}{\mathbf{F}}
\newcommand{\E}{\mathbb{E}}
\renewcommand{\S}{\mathbb{S}}

\newcommand{\R}{\mathbb{R}}

\newcommand{\I}{\mathcal{I}}
\newcommand{\II}{\mathrm{II}}
\newcommand{\M}{\mathbb{M}}
\newcommand{\Id}{\mathrm{Id}}
\renewcommand{\H}{\mathcal{H}}
\newcommand{\D}{\mathcal{D}}
\newcommand{\eps}{\epsilon}

\newcommand{\norm}[1]{\left\Vert#1\right\Vert}
\newcommand{\snorm}[1]{\Vert#1\Vert}
\newcommand{\abs}[1]{\left\vert#1\right\vert}
\newcommand{\set}[1]{\left\{#1\right\}}
\newcommand{\brac}[1]{\left(#1\right)}
\newcommand{\scalar}[1]{\left \langle #1 \right \rangle}

\DeclareMathOperator{\interior}{int}

\newenvironment{displayme}
  {\[ \begin{tabular}{m{0.95\textwidth} } \em
  \leftskip=0.5cm plus 0.5fil \rightskip=0.5cm plus -0.5fil \parfillskip=0cm plus 0.5\textwidth}
{\end{tabular}\]\ignorespacesafterend}

\begin{document}

\title{\Large Multi-Bubble Isoperimetric Problems}
\author{Emanuel Milman\thanks{Faculty of Mathematics, Technion-Israel Institute of Technology, Haifa 32000, Israel. Email: emilman@tx.technion.ac.il.}}

\begingroup
    \renewcommand{\thefootnote}{}
\endgroup

\date{}

\maketitle

\begin{abstract}
We survey recent advancements in the characterization of multi-bubble isoperimetric minimizers and the stability of soap bubble partitions. We conclude with some related open problems. 
\end{abstract}

\section{Modeling soap films}

The mathematical modeling of soap films has served as an impetus for the development of numerous facets of geometric measure theory. To a good approximation (the ``dry'' scenario), a soap film obtained by dipping a wire frame into soap water will locally minimize its surface area, yielding a minimal surface, as understood by Young, Laplace and Gauss in the first half of the 19th century. 
The first Fields medal was awarded in 1936 to J.~Douglas (jointly with L.~Ahlfors), who showed \cite{Douglas1931-PlateauProblem} (independently with T.~Rad\'o \cite{Rado1930-PlateauProblem}) the existence of a smoothly parametrized minimal surface whose boundary is a given Jordan curve (representing a wire frame) in $\R^n$, a problem dating back to Lagrange in 1760.
By construction, the Douglas--Rad\'o solutions are smoothly immersed topological discs, even though a soap film spanned by a simple smooth cycle may have higher genus or develop singularities, as observed in experiments by physicist J.~Plateau circa 1873. Nowadays, the Plateau problem entails establishing the existence and (partial) regularity of an $m$-dimensional minimal (generalized) surface spanning a given $(m-1)$-dimensional boundary in an $n$-dimensional ambient space. A solution is typically sought in an appropriate compact family of generalized surfaces which can accommodate various possible topologies, singularities, multiplicities, notions of boundary spanning, and other restrictions (such as orientation or the lack thereof).

We refer to \cite{MorganBook5Ed,David-ShouldWeSolvePlateau} and the references therein for a comprehensive discussion on various models for Plateau's problem in the context of soap films, as well as to \cite{KMS-PlateauAsLimitOfCapillarity,MaggiNovackRestrepo-PlateauBorders} for more recent developments. 
Fortunately, these subtle ambiguities in the mathematical modeling of soap films spanning a wire frame do not appear when modeling \emph{soap bubbles}, which are soap films enclosing trapped pockets of air (called bubbles). The number of bubbles $k$ is predetermined (and in our discussion, finite), and so are the $k$ volumes of trapped air. A stable configuration of $k$ soap bubbles is then a local minimizer of the total surface area of the soap film used to enclose the bubbles, given the $k$ volume constraints. This generalizes the single-bubble case $k=1$, where the round sphere is known since antiquity (at least in dimensions 2 and 3) to minimize surface area for a given volume by the classical isoperimetric inequality, and 
leads to a multi-bubble isoperimetric formulation, described next. 

\section{The multi-bubble isoperimetric setup}

A weighted Riemannian manifold $(M^n,g,\mu)$ consists of a smooth complete $n$-dimensional Riemannian manifold $(M^n,g)$ endowed with a measure $\mu$ having $C^\infty$ smooth positive density $\Psi$ with respect to the Riemannian volume measure $\vol_g$. The metric $g$ induces a geodesic distance on $M^n$, and the corresponding $k$-dimensional Hausdorff measure is denoted by $\H^k$. Let $\mu^k=e^{-W} \H^k$ and set the $\mu$-weighted volume to be $V_{\mu}:= \mu$. The $\mu$-weighted perimeter of a Borel subset $U \subset M$ of locally finite perimeter is defined as $A_\mu (U) : = \mu^{n-1} (\partial^* U)$, where $\partial^* U$ is the reduced (measure-theoretic) boundary of $U$ \cite{MaggiBook}. 

The Euclidean, spherical and hyperbolic model spaces $(M^n,g)$ are denoted by $\R^n$, $\S^n$ and $\HH^n$, respectively. They are endowed with their standard Riemannian volume measure $\mu = \vol_g$, and we will simply write $V$ and $A$ for volume and perimeter. Another important model space is the Gaussian one $\G^n$, obtained by endowing Euclidean space $\R^n$ with the standard Gaussian measure $\mu = \gamma^n := (2 \pi)^{-n/2} \exp(-|x|^2/2) dx$.

A $q$-partition $\Omega = (\Omega_1, \ldots, \Omega_q)$ of $(M,g,\mu)$ is a $q$-tuple of Borel subsets $\Omega_i \subset M$ having locally finite perimeter, such that $\{\Omega_i\}$ are pairwise disjoint and $V_\mu (M \setminus \cup_{i=1}^q \Omega_i) = 0$. Note that the sets $\Omega_i$, called cells, are not required to be connected. A $k$-tuple of pairwise disjoint cells $(\Omega_1,\ldots,\Omega_{k})$ so that $V_{\mu}(\Omega_i), A_{\mu}(\Omega_i) < \infty$ for all $i=1,\ldots,k$ is called a $k$-cluster, and its cells are called bubbles. Every $k$-cluster induces a partition by simply adding the ``exterior cell'' $\Omega_{k+1} := M \setminus \cup_{i=1}^{k} \Omega_i$; by abuse of notation, we will call the resulting $(k+1)$-partition  $\Omega = (\Omega_1,\ldots,\Omega_{k+1})$ a $k$-cluster (or $k$-bubble) as well. 
Of course, when $V_\mu(M) = \infty$ then necessarily $V_\mu(\Omega_{k+1}) = \infty$.

The $\mu$-weighted volumes vector $V_\mu(\Omega)$ and total perimeter (or surface area) $A_\mu(\Omega)$ of a $q$-partition $\Omega$ are defined as
\begin{align*}
	V_{\mu} (\Omega) & := ( V_{\mu} (\Omega_1), \ldots, V_{\mu} (\Omega_q) )  , \\ 	A_{\mu} (\Omega) & := \frac{1}{2} \sum_{i=1}^q A_{\mu} (\Omega_i) = \sum_{1\leq i <j \leq  q} \mu^{n-1} (\Sigma_{ij}),
\end{align*}
where
\[
 \Sigma_{ij}:= \partial^* \Omega_i \cap \partial^* \Omega_j
 \]
 denotes the $(n-1)$-dimensional  interface between cells $\Omega_i$ and $\Omega_j$.  We set $\Delta^{\hspace{-1pt}(q-1)}[T] := \{ v \in [0,\infty)^{q-1} \times [0,\infty] : \sum_{i=1}^q v_i =T \}$, where $T = V_\mu(M)$.

The isoperimetric problem for $k$-clusters consists of identifying those clusters $\Omega$ of prescribed volume $V_\mu (\Omega) = v \in \interior \Delta^{\hspace{-1pt}(k)}[T]$ which minimize the total perimeter $A_{\mu}(\Omega)$; local minimizers are also interesting to classify. By modifying an isoperimetric minimizing cluster on a null set, we may and will assume that its cells $\Omega_i$ are open and satisfy $\overline{\partial^* \Omega_i} = \partial \Omega_i$. 

The solutions to the classical isoperimetric problem, corresponding to the single-bubble case $k=1$, constitute some of the most beautiful and ancient results in geometry, and play a key role in various facets of differential geometry, analysis, PDE, calculus of variations, geometric measure theory, probability, mathematical physics, boolean analysis, combinatorics, etc...~It is well known that geodesic balls $\Omega_1$ of prescribed volume uniquely minimize perimeter on all model spaces $\R^n$, $\S^n$ and $\HH^n$ \cite{Steiner-1842,Schwarz-1890,SchmidtIsopOnModelSpaces,BuragoZalgallerBook}. It is also classical that halfplanes $\Omega_1$ of prescribed Gaussian volume uniquely minimize Gaussian perimeter on $\G^n$ \cite{SudakovTsirelson,Borell-GaussianIsoperimetry,CarlenKerceEqualityInGaussianIsop}; their flat boundaries can be thought of as degenerate flat spheres. Consequently, we will collectively refer to complete constant curvature hypersurfaces in $\M^n \in \{\R^n, \S^n, \HH^n\}$ as ``generalized spheres'' -- in $\R^n$ these are spheres and hyperplanes, and in $\HH^n$ these are geodesic spheres, horospheres, and equidistant hypersurfaces. 

The multi-bubble isoperimetric problem for $k$-clusters (when $k \geq 2$) already poses a much greater challenge. Even just formulating a reasonable conjecture for general $k$ requires some ingenuity.

\section{Isoperimetric conjectures}

When the number of bubbles $k$ is much larger than the ambient dimension $n$, it is entirely unclear what could be a plausible minimizing configuration in $\R^n$, $\S^n$, $\HH^n$ or $\G^n$. Even in $\R^3$, computer simulations using K.~Brakke's surface evolver suggest that a minimizing $6$-bubble may not have spherical interfaces \cite{SullivanOldSurvey}. A more interesting and tractable question would be, say in the equal-volume case in $\R^2$, to ask what is the asymptotic behaviour when $k \rightarrow \infty$, or what is an optimal tiling (see e.g. \cite{CoxGranerEtAl,CoxMorganGraner,Hales-Honeycomb, HeppesMorgan,CarocciaMaggi-StableHales,CarocciaDeMasonMaggi-AlmostHoneycomb, PeriodicDoubleTilings}), but we do not expand on this here. Instead, let us recall the following definition and corresponding conjectures in the case that $k \leq n+1$, which were put forth by J.~Sullivan in the 1990's \cite[Problem 2]{OpenProblemsInSoapBubbles96}. 
\medskip

A centered simplicial $q$-partition $\Omega^{q}_0[\R^N]$ of $\R^N$, $2 \leq q \leq N+1$, is given by the Voronoi cells
\[
\Omega^{q}_{0,i}[\R^N] := \set{p \in \R^N : \argmin_{j=1,\ldots,q} \scalar{p,\c_j} = \{i\} } ~,~ i = 1,\ldots,q ,
\]
where $\{\c_j\}_{j=1,\ldots,q}$ are $q$ (distinct) equidistant points in $\R^{N}$ with $\sum_{j=1}^q \c_j = 0$. Of course, 
one cannot find more than $N+1$ equidistant points in $\R^{N}$.

\begin{figure}[htbp]
\centering
\includegraphics[scale=0.52]{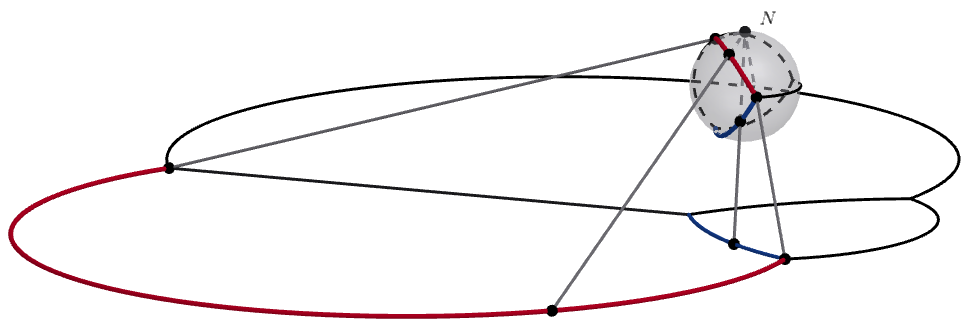}
\caption{ \label{fig:stereographic}
Stereographic projection of an equal-volume partition of $\S^2$ into 4 Voronoi cells, yielding a standard triple-bubble in $\R^2$. 
}
\end{figure}

\begin{figure}
    \begin{center}
            \raisebox{-0.1\height}{\includegraphics[scale=0.1]{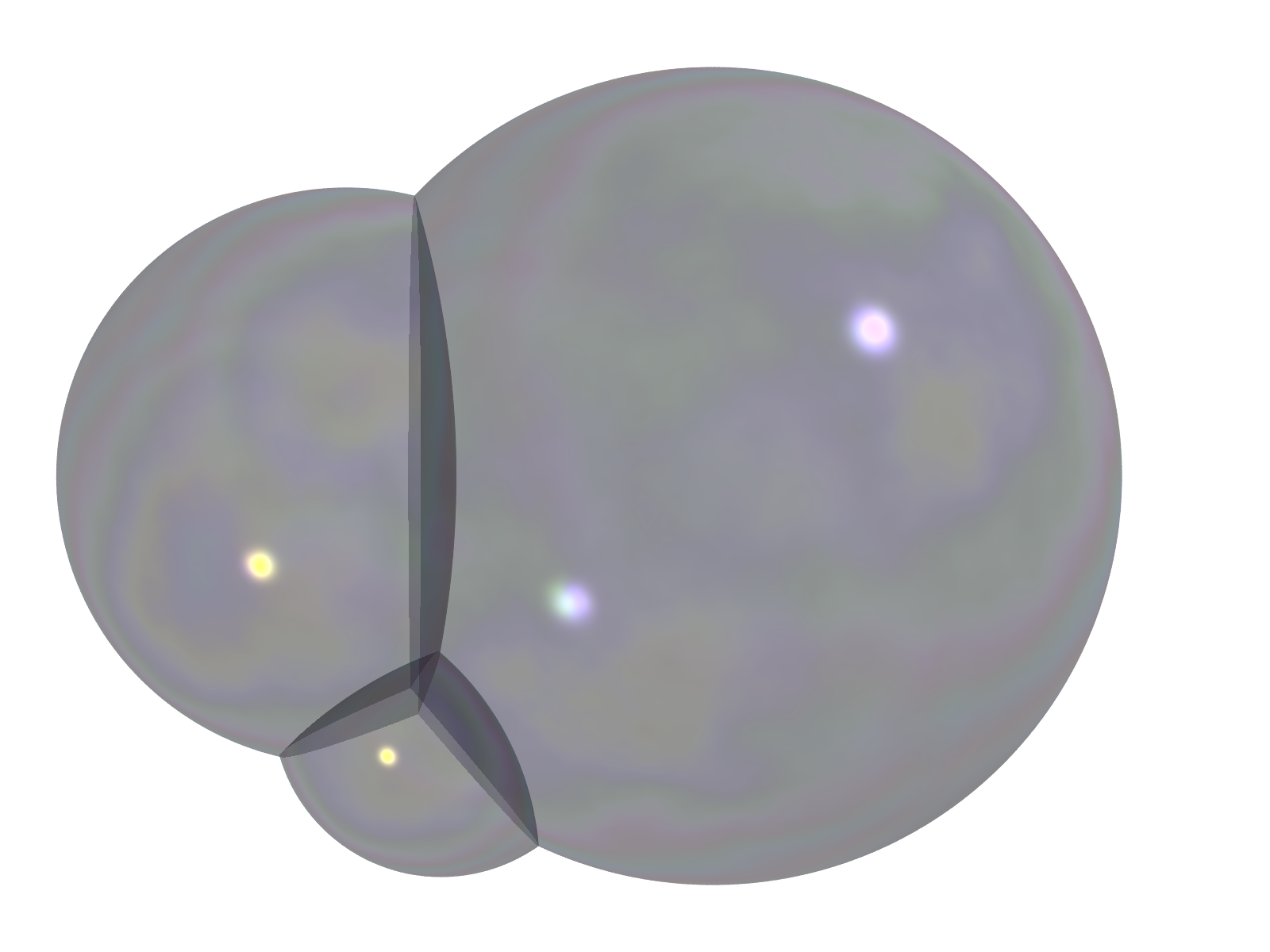}}
        \hspace{20pt}
        \begin{tikzpicture}[scale=1.2]
            \input{triple-bubble}
        \end{tikzpicture}
     \end{center}
     \caption{
         \label{fig:triple-bubble}
        Left: a standard triple-bubble in $\R^3$. Right: the $2$D cross-section through its plane of symmetry. 
     }
\end{figure}

\begin{figure}
    \begin{center}
        \includegraphics[scale=0.13]{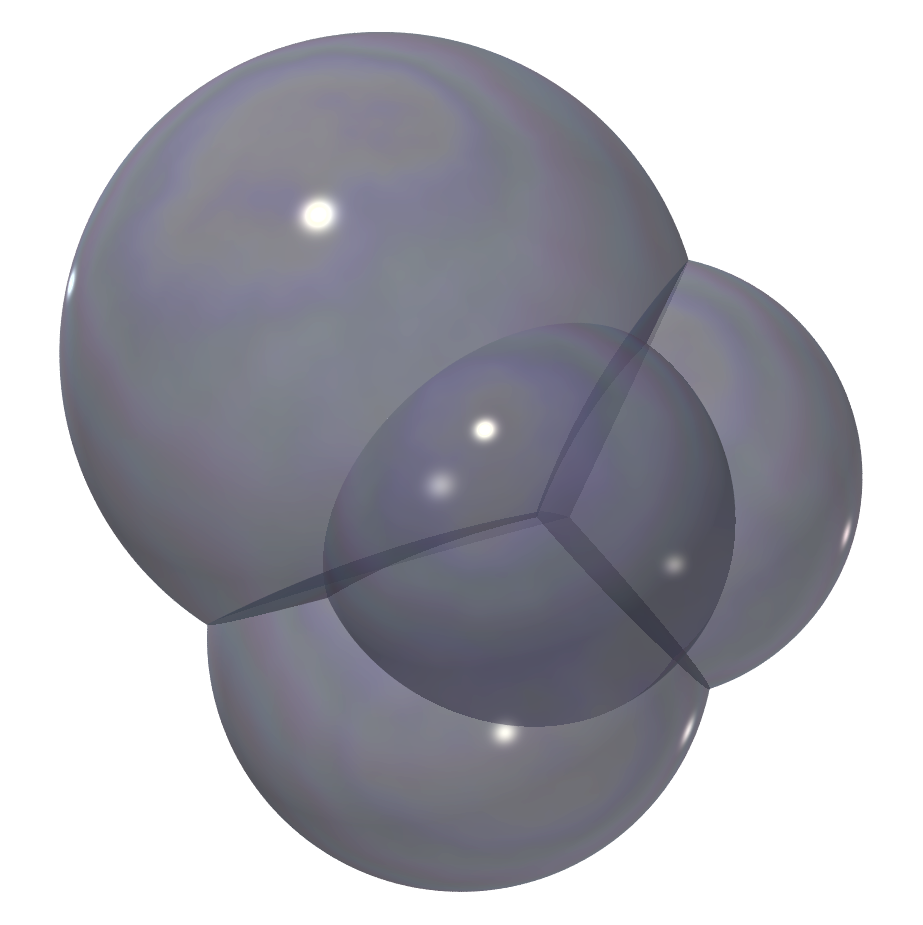}
        \includegraphics[scale=0.13]{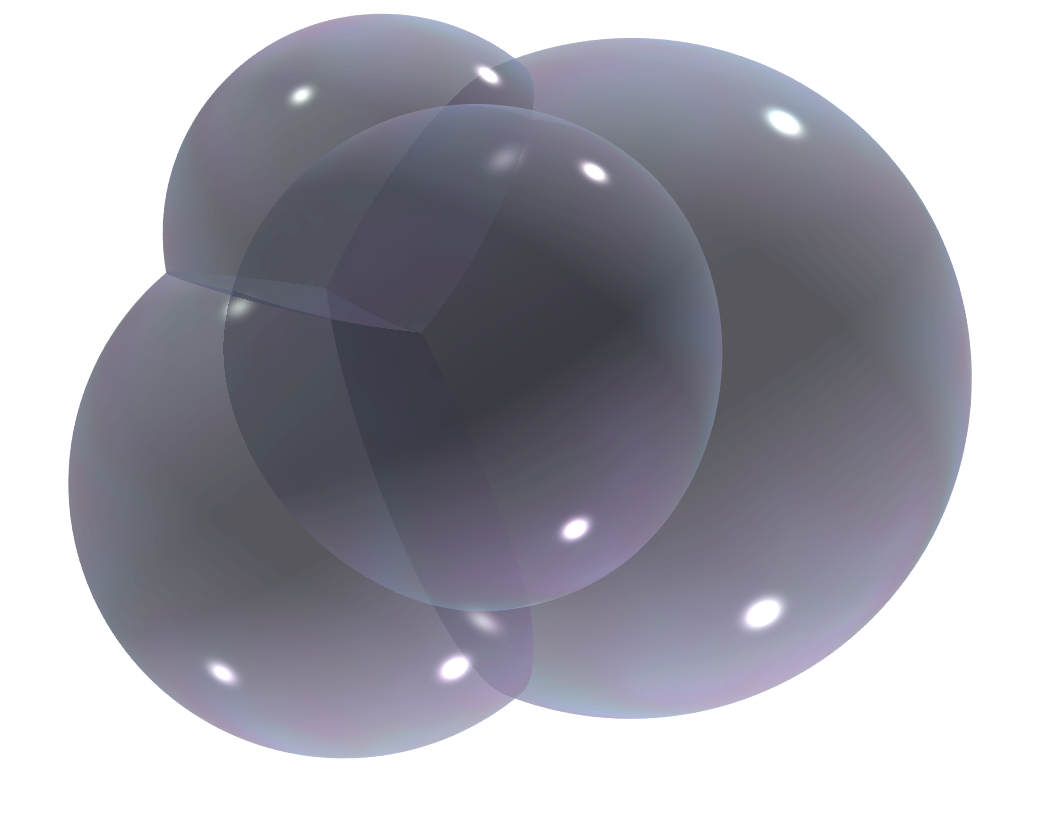}
        \includegraphics[scale=0.13]{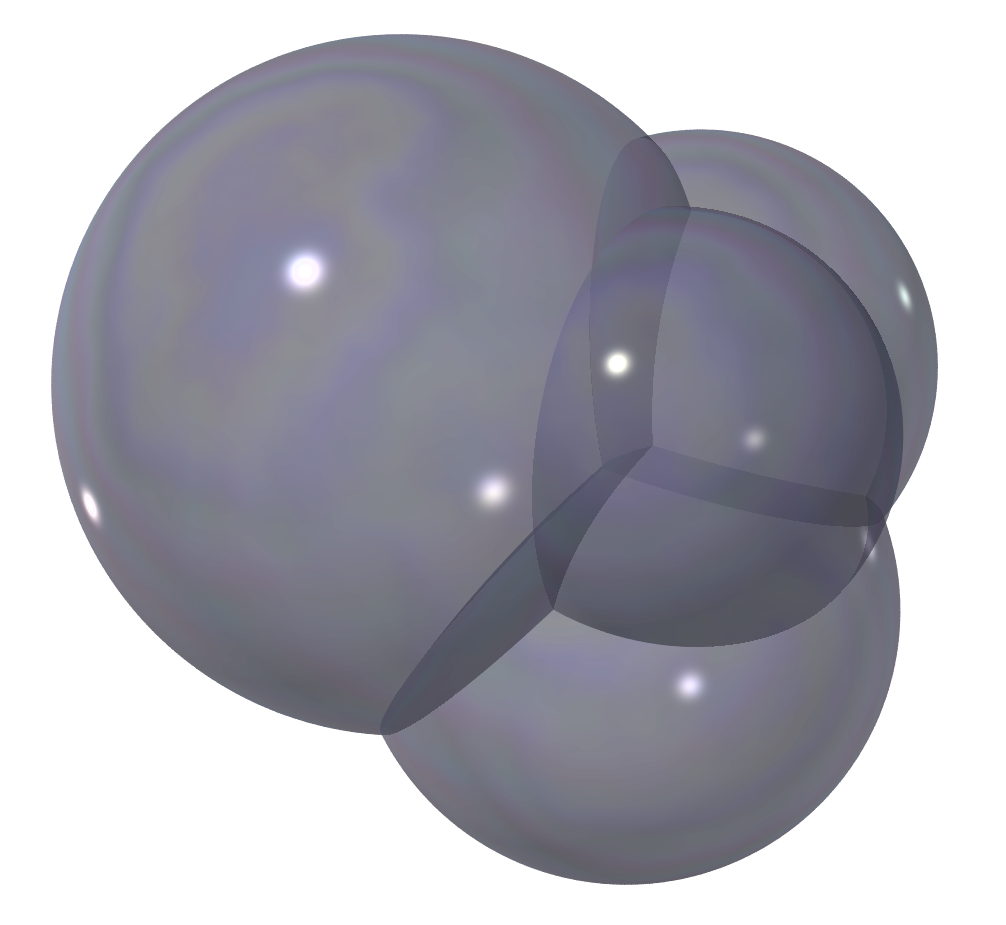}
             \end{center}
     \caption{
         \label{fig:quadruple-bubble}
        A standard quadruple-bubble in $\R^3$ (also, the cross-section of a standard quadruple-bubble in $\R^4$ through its hyperplane of symmetry) from different angles. 
     }
\end{figure}

\begin{definition}[Standard partitions and bubbles in $\R^n$, $\S^n$ and $\HH^n$]
Given $1 \leq k \leq n+1$, an equal-volume standard $(k+1)$-partition $\Omega^{k+1}_0[\S^n]$ of $\S^n$  is given by
\[
\Omega^{k+1}_{0,i}[\S^n] := \Omega^{k+1}_{0,i}[\R^{n+1}] \cap \S^n ~,~ i = 1,\ldots,k+1 . 
\]
A standard $(k+1)$-partition $\Omega^{k+1}[\R^n]$ of $\R^n$ is a stereographic projection of $\Omega^{k+1}_0[\S^n]$. A standard $(k+1)$-partition $\Omega^{k+1}[\M^n]$ of $\M^n \in \{\S^n,\HH^n\}$ is a stereographic projection of a corresponding one in $\R^n$. The resulting partition is a standard $k$-bubble in $\M^n \in \{\R^n,\S^n,\HH^n\}$ if $V(\Omega^{k+1}_i[\M^n]) \in (0,\infty)$ for all $i=1,\ldots,k$. See Figures \ref{fig:stereographic}, \ref{fig:triple-bubble} and \ref{fig:quadruple-bubble}. 
\end{definition}

\begin{figure}
\begin{center}
\begin{tikzpicture}
    \coordinate(O) at (0, 0.0);   
    \coordinate(A) at ({-1.3 * sqrt(3)}, 1.3);
    \coordinate(B) at ({1.3 * sqrt(3)}, 1.3);
    \coordinate(C) at (0, -2.1);   
    
    \draw[thick] (O) -- (A); 
    \draw[thick] (O) -- (B); 
    \draw[thick] (O) -- (C); 
\end{tikzpicture}
\hspace{50pt}
\includegraphics[scale=0.27]{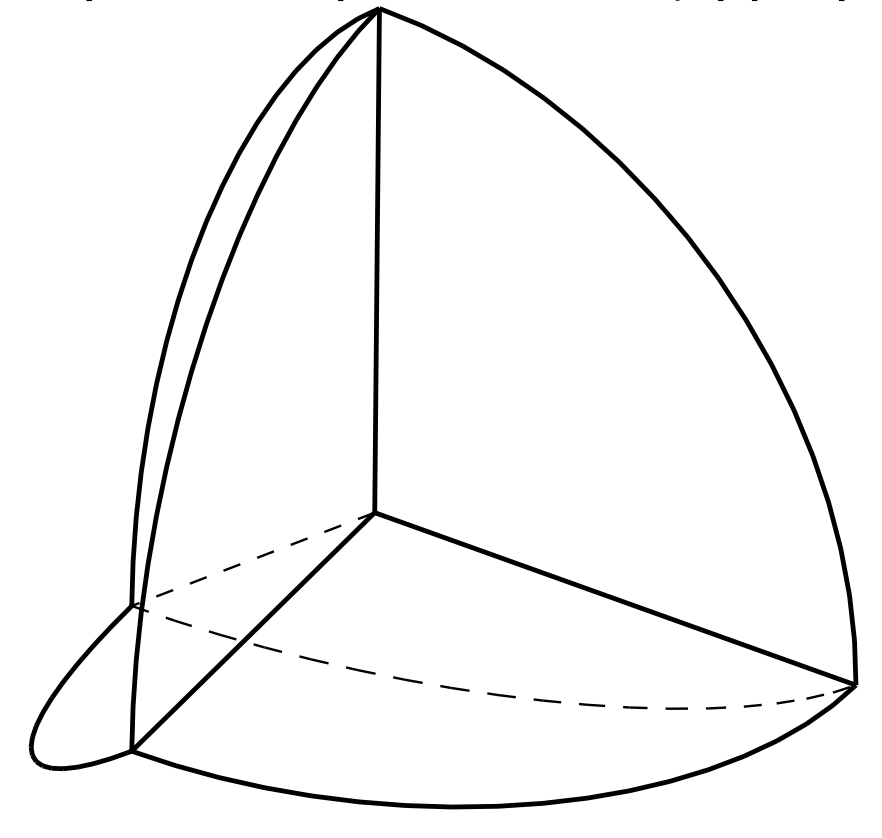}
\end{center}
\caption{
         \label{fig:simplicial-bubble}
        A standard simplicial double-bubble in $\G^2$, whose boundary is a $\Y$-cone  (left), and a standard simplicial triple-bubble in $\G^3$, whose boundary is a $\T$-cone (right). 
     }
\end{figure}

Note that all standard partitions are generated from the equal volumes one by the group of M\"obius transformations. Indeed, the composition of two stereographic projections back and forth from (the one point compactification of) $\R^n$ is a M\"obius transformation, obtained by composing isometries with scaling and spherical inversion. By Liouville's classical theorem, all such global conformal automorphisms of $\R^n \cup \{\infty\}$ and $\S^n$ are given by M\"obius transformations when $n \geq 3$. However, on $\HH^n$, there are no non-trivial global conformal automorphisms (besides isometries), and so a bit more care is required when constructing standard bubbles, and a stereographic projection from a hemisphere $\S^n_+$ is actually used (see \cite{EMilmanXu-Stability}). It is well known that the group of all M\"obius automorphisms (say on $\S^n$) is isomorphic to the Lorentz group $SO_1(n + 2)$ of (orientation preserving) isometries of Minkowski space-time $\R^{n+2}$ with signature $(+,\ldots,+,-)$, making a link to projective geometry.

 Since stereographic projection preserves generalized spheres, the interfaces of all standard partitions are indeed subsets of generalized spheres (as for a single-bubble minimizer). In addition, stereographic projection is a conformal diffeomorphism, preserving both angles and the combinatorial incidence structure, and hence all $m$-tuples of cells for $m \leq n+1$ meet at equal angles like the cones over the $m$ facets of a regular $(m-1)$-dimensional simplex. In particular, every pair of cells $\Omega_i$ and $\Omega_j$ meet along a non-empty interface $\Sigma_{ij} \neq \emptyset$. Furthermore, the $3$ interfaces $\Sigma_{ab}$, $\Sigma_{bc}$ and $\Sigma_{ca}$ meet at $120^{\circ}$ angles (a $\Y$ singularity), and the $6$ interfaces $\{\Sigma_{ij}\}_{\{i,j\} \subset \{a,b,c,d\}}$ meet like the cone over the edges  of a regular tetrahedron (a $\T$ singularity) -- the only two types of singularities observed experimentally in $\R^3$ by Plateau, referred to as ``Plateau's laws''.

\begin{conjecture}[Multi-Bubble Isoperimetric Conjecture on $\M^n \in \{\R^n, \S^n, \HH^n\}$]
For all $2 \leq k \leq n+1$, a standard $k$-bubble uniquely minimizes total perimeter among all $k$-clusters $\Omega$ on $\M^n$ of prescribed volume $V(\Omega) = v \in \interior \Delta^{\hspace{-1pt}(k)}[V(\M^n)]$. 
\end{conjecture}

The natural analogue in $\G^n$ of the above construction is to use a standard \emph{simplicial} bubble, whose interfaces are flat (as for a single-bubble Gaussian minimizer) \cite{CorneliCorwinEtAl-DoubleBubbleInSandG, Schechtman-ApproxGaussianMultiBubble}.

\begin{definition}[Standard simplicial bubbles in $\G^n$]
Given $1 \leq k \leq n$, a standard simplicial $k$-bubble $\Omega^{k+1}[\G^n]$ in $\G^n$ is a $(k+1)$-partition obtained by translating a centered simplicial partition $\Omega^{k+1}_0[\R^n]$. See Figure \ref{fig:simplicial-bubble}. 
\end{definition}

\begin{conjecture}[Multi-Bubble Isoperimetric Conjecture on $\G^n$]
For all $2 \leq k \leq n$, a standard simplicial $k$-bubble uniquely minimizes total Gaussian perimeter among all $k$-clusters $\Omega$ on $\G^n$ of prescribed Gaussian volume $V_{\gamma^n}(\Omega) = v \in \interior \Delta^{\hspace{-1pt}(k)}[1]$. 
\end{conjecture}

That a standard $k$-bubble in $\R^n$ of prescribed volume $v$ exists and is unique (up to isometries) for all $v \in \interior \Delta^{\hspace{-1pt}(k)}[\infty]$ and $2 \leq k \leq n+1$ was shown by Montesinos-Amilibia \cite{MontesinosStandardBubbleE!} (see \cite{EMilmanNeeman-TripleAndQuadruple,EMilmanNeeman-GaussianMultiBubble} for the corresponding analogous statements on $\S^n$ and $\G^n$). 

\section{Partial answers} \label{sec:partial}

When $n=2$, the double-bubble conjectures (case $k=2$) on the various model spaces described above are well understood and fully resolved. On $\R^2$, $\S^2$ and $\HH^2$ this was established in \cite{SMALL93}, \cite{Masters-DoubleBubbleInS2} and \cite{CottonFreeman-DoubleBubbleInSandH}, respectively; the case of $\G^2$ was only recently resolved (see below). The triple-bubble conjectures (case $k=3$) on $\R^2$ and $\S^2$ were established by Wichiramala \cite{Wichiramala-TripleBubbleInR2} and Lawlor \cite{Lawlor-TripleBubbleInR2AndS2}, respectively, but to the best of our knowledge has not been worked out on $\HH^2$ and remains open. While this falls outside the scope of Sullivan's conjecture, we also mention the work of Paolini, Tamagnini and Tortorelli \cite{PaoliniTamagnini-PlanarQuadraupleBubbleEqualAreas,PaoliniTortorelli-PlanarQuadrupleEqualAreas}, who have identified a unique minimizing quadruple-bubble of equal volumes $v_1=v_2=v_3=v_4$ in the plane $\R^2$.

In a landmark work, the double-bubble conjecture on $\R^3$ was confirmed by Hutchings, Morgan, Ritor\'e and Ros \cite{DoubleBubbleInR3-Announcement,DoubleBubbleInR3} following prior contributions in \cite{HHS95,HassSchlafly-EqualDoubleBubbles,Hutchings-StructureOfDoubleBubbles}, and was subsequently extended to all $\R^n$ \cite{SMALL03,Reichardt-DoubleBubbleInRn,Lawlor-DoubleBubbleInRn}. On $\S^n$ for $n\geq 3$ and $\G^n$ for $n \geq 2$, only partial results for the double-bubble problem were known until recently \cite{CottonFreeman-DoubleBubbleInSandH, CorneliHoffmanEtAl-DoubleBubbleIn3D,CorneliCorwinEtAl-DoubleBubbleInSandG}.

\smallskip

In \cite{EMilmanNeeman-GaussianMultiBubble}, in collaboration with Joe Neeman, we fully resolved the multi-bubble conjecture on $\G^n$ in the Gaussian setting (including the case when $k=n$, the largest value to which the conjecture applies):
\begin{displayme}
The multi-bubble isoperimetric conjecture on $\G^n$ holds true for the entire applicable range $2 \leq k \leq n$. 
\end{displayme}
\noindent
Subsequently, together with Neeman in \cite{EMilmanNeeman-TripleAndQuadruple,EMilmanNeeman-QuintupleBubble}, we were able to establish the following:
\begin{displayme}
The double-, triple-, quadruple- and quintuple-bubble conjectures (cases $k=2,3,4,5$) hold true on $\M^n$ for $\M^n \in \{\R^n,\S^n\}$ and $n \geq k$ (without uniqueness on $\R^n$ in the quintuple case). 
\end{displayme}

The remaining cases on $\R^n$, $\S^n$ and $\HH^n$ are still open. In particular, unlike in the Gaussian setting, we are unable to handle the largest value $k = n+1$ to which the conjectures apply even when $k \leq 5$ (such as the quadruple-bubble case $k=4$ in dimension $n=3$). One exception is the equal-volumes case $v_1 = \ldots = v_{k+1}$ of the multi-bubble conjecture on $\S^n$ for all $2 \leq k \leq n+1$, which follows easily from the equal-volumes case on $\G^{n+1}$. Indeed, by rotation invariance of the Gaussian measure, both (weighted) volume and perimeter on $\S^n$ and $\G^{n+1}$ coincide for \emph{centered} cones (after normalizing $V(\S^n)$ to be $1$), and the unique equal-volumes minimizer in $\G^{n+1}$ for all $1 \leq k \leq n+1$ is the centered standard simplicial bubble (whose cells are centered cones).

Note that the hyperbolic case $\HH^n$ is presently out of reach of our methods. 
However, we are able to obtain some additional partial results on $\R^n$ and $\S^n$.

\begin{definition}[Spherical Voronoi partition of $\S^n$] 
A $q$-partition $\Omega$ of $\S^n$, all of whose cells are non-empty, 
is called a spherical Voronoi partition if there exist $\{ \c_i \}_{i=1,\ldots,q} \subset \R^{n+1}$ and $\{ \k_i \}_{i=1,\ldots,q} \subset \R$ so that the following holds:
\begin{enumerate}
\item
For every non-empty interface $\Sigma_{ij} \neq \emptyset$, $\Sigma_{ij}$ is a relatively open subset of a geodesic sphere $S_{ij}$ in $\S^n$ with quasi-center $\c_{ij} = \c_i - \c_j$ and curvature $\k_{ij} = \k_i - \k_j$. \\
The quasi-center $\c$ of a geodesic sphere $S$ is the vector $\c := \n - \k p$ at any of its points $p \in S$ (where $\k$ is the curvature with respect to the unit normal $\n$). 
\item 
The following Voronoi representation holds:
\[ 
\Omega_i  = \set{ p \in \S^n : \argmin_{j=1,\ldots,q} \scalar{p,\c_j} + \k_j = \{ i\}  } =  \bigcap_{j \neq i} \; \set{ p \in \S^n : \scalar{p,\c_{ij}} + \k_{ij} <  0 } .
\] 
\end{enumerate}
\end{definition}

\begin{definition}[Spherical Voronoi partition of $\M^n \in \{ \R^n, \HH^n \}$]
A $q$-partition $\Omega$ of $\M^n$ is called spherical Voronoi, if it is a stereographic projection of a spherical Voronoi $q$-partition of $\S^n$. 
\end{definition}

\begin{figure}
    \begin{center}
         \hspace*{37pt}
        \includegraphics[scale=0.108]{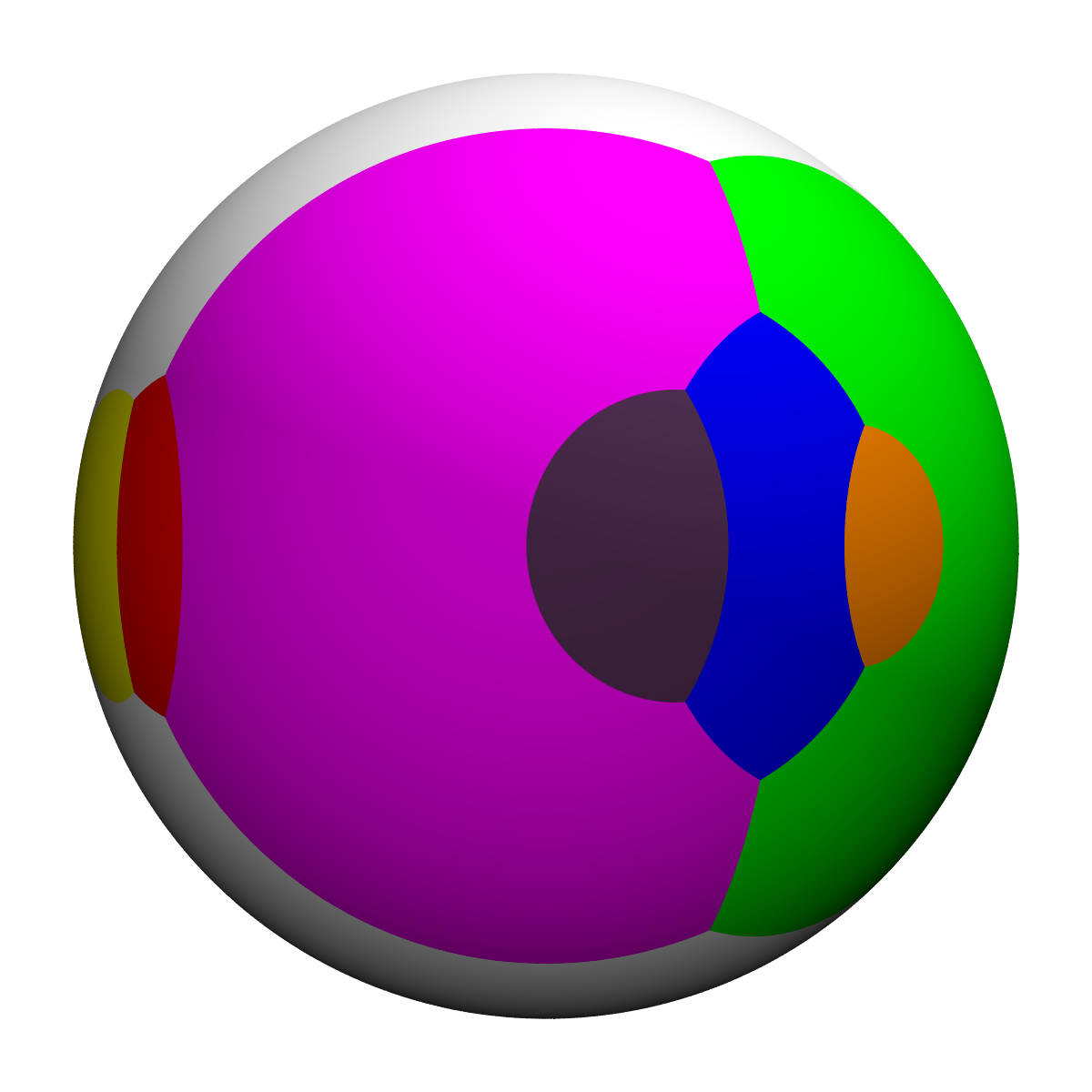}
        \hspace{30pt}
        \raisebox{0.01\height}{
        \begin{tikzpicture}[scale=1.76]
            \input{cluster-stereographic-color}
        \end{tikzpicture}
        }
        \newline
                \includegraphics[scale=0.108]{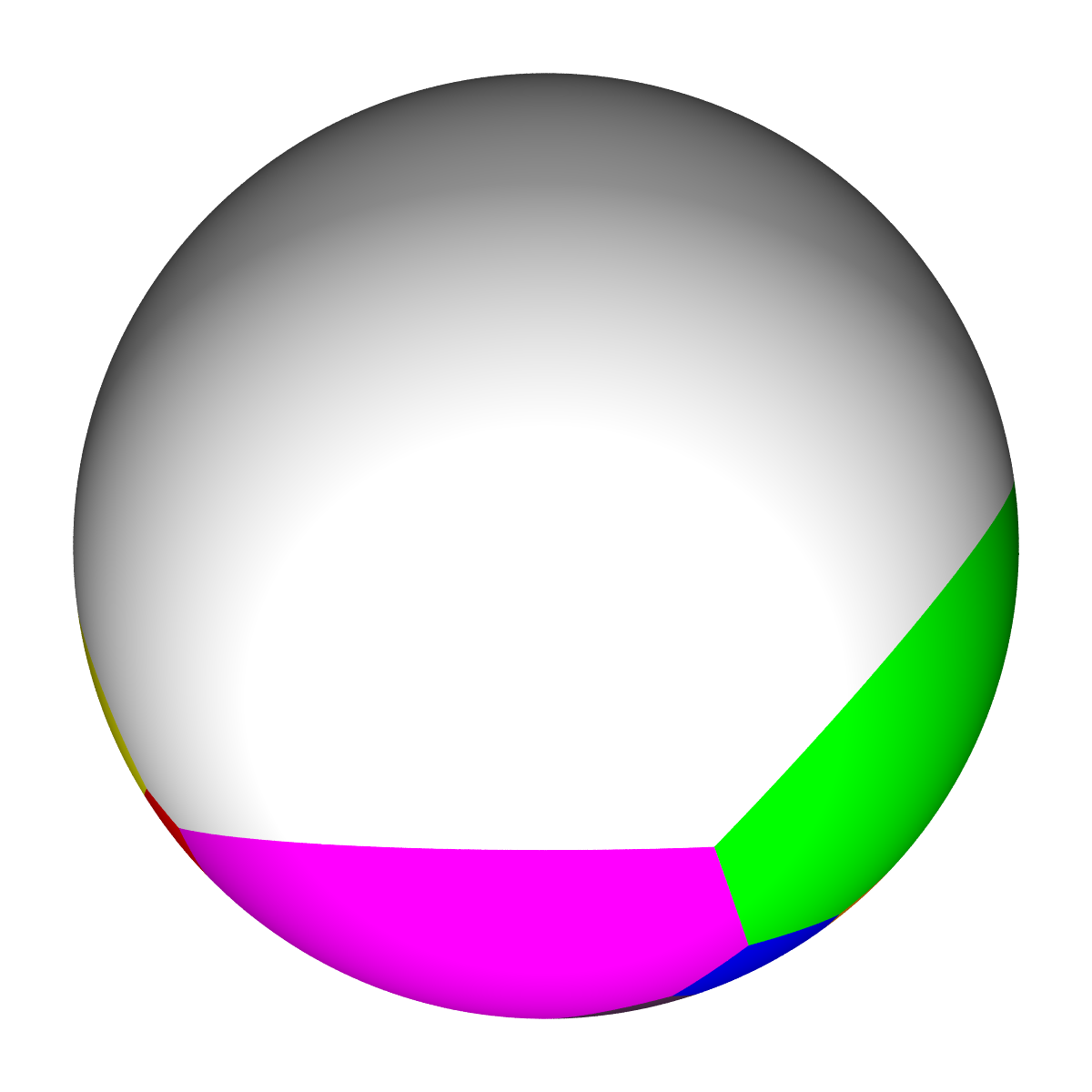}
        \hspace{60pt}
                \raisebox{0.08\height}{
        \begin{tikzpicture}[scale=1.97]
            \input{cluster-brep-color}
        \end{tikzpicture}
        }
     \end{center}
     \caption{
         \label{fig:SphericalVoronoi}
         Top left: A spherical Voronoi cluster $\Omega^\S$ in $\S^2$; Top right: A spherical Voronoi cluster $\Omega^\R$ in $\R^2$ obtained from $\Omega^\S$ by stereographic projection; Bottom left: $\Omega^\S$ drawn from above; Bottom right: the orthogonal projection of $\Omega^\S$ onto its plane of symmetry consists of convex polyhedral cells (colors lightened for better contrast). 
     }
\end{figure}

\begin{remark} \label{rem:intro-standard-bubbles}
Spherical Voronoi partitions are closed under M\"obius transformations. Consequently, all standard $k$-bubbles in $\R^n$, $\S^n$ and $\HH^n$ are spherical Voronoi clusters. Note that each cell of a spherical Voronoi partition of $\S^n$ is the intersection of an open convex polyhedron in $\R^{n+1}$ with $\S^n$. See Figure \ref{fig:SphericalVoronoi}. 
\end{remark}

\begin{theorem}[Spherical Voronoi structure and connectedness of cells \cite{EMilmanNeeman-TripleAndQuadruple}] \label{thm:Voronoi}
Let $\Omega$ be an isoperimetric minimizing $k$-cluster in $\M^n \in \{\R^n, \S^n\}$, $k \leq n$, with $V(\Omega)\in \interior  \Delta^{\hspace{-1pt}(k)}[V(\M^n)]$. Then $\Omega$ is a spherical Voronoi cluster, and all of its cells are connected. 
\end{theorem}
 
In particular, this resolved a conjecture of A.~Heppes \cite[Problem 5]{OpenProblemsInSoapBubbles96} on the connectedness of the cells, and the question of whether there can be empty chambers trapped by minimizing bubbles \cite[Chapter 13]{MorganBook5Ed}. Showing the connectedness of the cells was in fact the main difficulty in the resolution of the double-bubble conjecture on $\R^n$ in \cite{DoubleBubbleInR3,SMALL03,Reichardt-DoubleBubbleInRn}.  
While it is not clear whether splitting a bubble into several connected components of air is physically possible, such a splitting might \emph{a priori} lead to savings in the total perimeter of the configuration, and so this needs to be \emph{a posteriori} ruled out.

The spherical Voronoi structure of a minimizing $k$-cluster reduces the multi-bubble isoperimetric problem (on $\M^n \in \{\R^n, \S^n\}$, when $k \leq n$) to a finite-dimensional configuration space, governed by the quasi-center $\{ \c_i \}_{i=1,\ldots,k+1}$ and curvature $\{ \k_i \}_{i=1,\ldots,k+1}$ parameters. In fact, for all $k \leq n+1$, standard $k$-bubbles in $\M^n$ are characterized as those spherical Voronoi $k$-clusters for which the interfaces $\Sigma_{ij}$ are non-empty for all $1 \leq i < j \leq k+1$ \cite{EMilmanNeeman-TripleAndQuadruple}. Consequently, in order to resolve the multi-bubble isoperimetric conjectures on $\M^n \in \{\R^n, \S^n\}$ when $k \leq n$, it remains to show that the cell-incidence graph of a minimizing cluster must be the complete graph, thereby reducing the problem to a combinatorial one.

\section{The local isoperimetric problem}

When $V_{\mu}(M^n)=\infty$, such as for $M^n \in \{\R^n,\HH^n \}$, one can also consider the \emph{local} isoperimetric problem for $q$-partitions $\Omega$ with prescribed volume $V_\mu(\Omega) = v$ when at least two of the prescribed volumes are infinite ($\exists i \neq j$ with $v_i=v_j = \infty$). In that case, $A_\mu(\Omega) = \infty$ and so the global minimization problem does not make sense. Instead, one considers \emph{locally} minimizing partitions, which minimize the total relative perimeter in any bounded open $K \subset M^n$ among all competing $q$-partitions $\Omega'$ with $V_\mu(\Omega') = V_\mu(\Omega)$ so that $\Omega'_i \Delta \Omega_i \Subset K$ for all $i$. 
This line of investigation was pioneered by Alama, Bronsard and Vriend \cite{ABV-LensClusters} on $\R^2$ for the case $v = (1,\infty,\infty)$, and systematically studied by Novaga, Paolini and Tortorelli in \cite{NPT-LocallyIsoperimetricPartitions}. Using a limit argument and a closure theorem, it was shown in \cite[Section 3]{NPT-LocallyIsoperimetricPartitions} that the results 
from the preceding section imply that standard $q$-partitions of $\R^n$ for $q=2,3,4$ ($n \geq 2$), $q=5$ ($n \geq 4$) and $q=6$ ($n \geq 5$) are locally minimizing (but does not exclude the existence of other non-standard partitions which are also locally minimizing when $n \geq 3$). 
Uniqueness of the local minimizing $q$-partitions of $\R^2$ for $q \leq 4$ was established in \cite{NPT-LocallyIsoperimetricPartitions}, and of the case $v=(1,\infty,\infty)$ on $\R^n$ for all $n \leq 7$ in \cite{BronsardNovak-DifferentTensions}, but the latter uniqueness fails to hold for $8 \leq n \leq 2700$ \cite{BronsardNovackEtAl-NonUniquenessOfLens,NPT-NonUniquenessOfLensInR8}. 

Some additional isoperimetric results for clusters and partitions are obtained in 
\cite{CLM-StableDoubleBubbleInR2, PaoliniTamagnini-PlanarQuadraupleBubbleEqualAreas,PaoliniTortorelli-PlanarQuadrupleEqualAreas,NPST-ClustersViaConcentrationCompactness, NPST-ClustersWithInfinitelyManyCells,DeRosaTione-ConvexStationaryBubbles,BronsardNovak-DifferentTensions,BCT-StableLensPartitionInR2,PeriodicDoubleTilings}.

\section{Stability} \label{sec:stability}

Showing the global (or local) minimality of the conjectured clusters (or partitions) in the remaining range of parameters seems to be a difficult task, and so confirmation of their local minimality, even in an infinitesimal sense, is already quite challenging, and provides valuable evidence for the validity of the conjectures. 
A standard way in the calculus of variations to probe the infinitesimal minimality of a given configuration is to test the non-negativity of the second variation (modulo the volume constraint), a property called ``stability'' (see Sections \ref{sec:GMT} and \ref{sec:proofs} for a more concrete definition). 

\begin{figure}
    \begin{center}
  \includegraphics[scale=0.28]{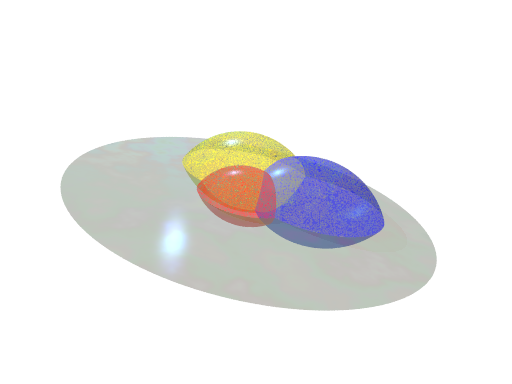}
     \includegraphics[scale=0.28]{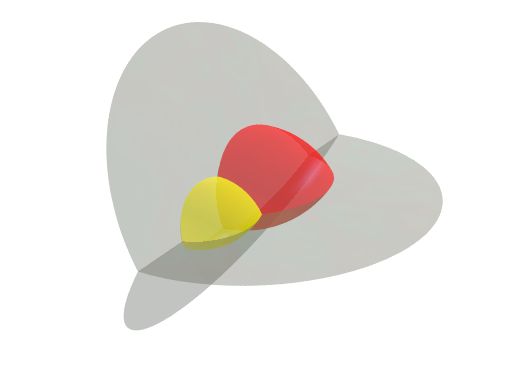}
     \includegraphics[scale=0.28]{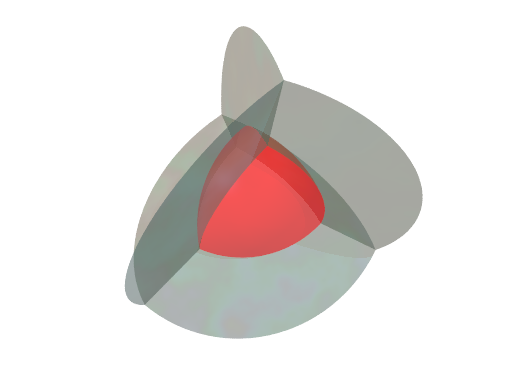}
     \end{center}
     \vspace{-10pt}
     \caption{
         \label{fig:standard-partitions}
         Standard $5$-partitions of $\R^3$ with $2$ (left), $3$ (middle) and $4$ (right) unbounded cells. All three are conjectured to be locally minimizing perimeter under volume constraint, but they are only known to be stable. 
     }
\end{figure}

Modulo technicalities, we confirmed in \cite{EMilmanXu-Stability}, in collaboration with Botong Xu, that:
\begin{displayme}
For all $1 \leq k \leq n+1$ and $n \geq 3$, standard $k$-bubbles in $\R^n$, $\S^n$ and $\HH^n$ are stable.
\end{displayme}
Moreover, this actually holds for all standard partitions, such as the ones depicted in Figure \ref{fig:standard-partitions}:
\begin{displayme}
For all $2 \leq q \leq n+2$ and $n \geq 3$, standard $q$-partitions of $\R^n$, $\S^n$ and $\HH^n$ are stable.
\end{displayme}
Our stability results are not restricted to standard partitions. A partition is called flat if all of its interfaces $\Sigma_{ij}$ are flat (i.e.~totally geodesic). A partition $\Omega$ of $\S^n$ is called \emph{M\"obius-flat} if  there exists a M\"obius automorphism $T : \S^n \rightarrow \S^n$ so that $T\Omega$ is flat. 
The same terminology is used on $\R^n$ and $\HH^n$ if this holds after a stereographic projection to $\S^n$ or the hemisphere $\S^n_+$, respectively. 
In \cite{EMilmanXu-Stability}, we showed that:
\begin{displayme}
For all $q \geq 2$ and $n \geq 3$, regular, M\"obius-flat, spherical Voronoi, $q$-partitions of $\R^n$, $\S^n$ and $\HH^n$ are stable.
\end{displayme}
We defer the precise definition of \emph{regularity} to Section \ref{sec:GMT}, but for now only mention that a regular partition combinatorially obeys Plateau's laws -- three cells meet like the cone over the vertices of a triangle, and four cells meet like the cone over the edges of a tetrahedron.

\begin{figure}
    \begin{center}
         \includegraphics[scale=0.3]{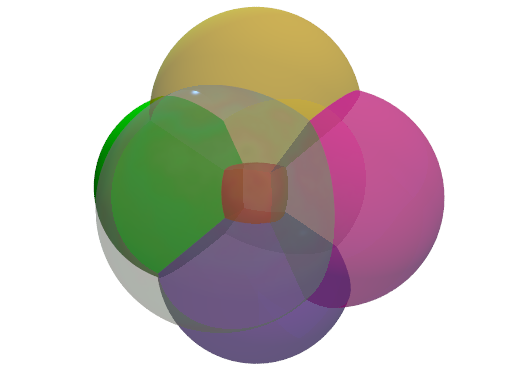}
            \includegraphics[scale=0.3]{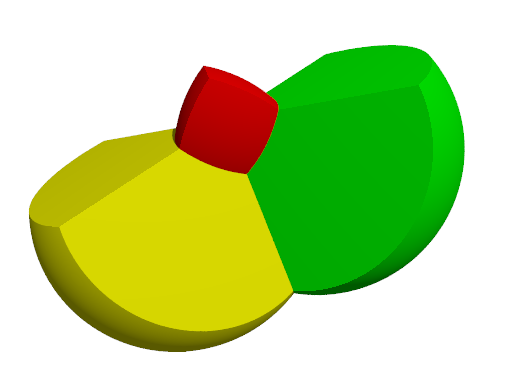}
     \end{center}
     \vspace{-10pt}
     \caption{
         \label{fig:nonstandard-partition1}
         A non-standard 7-bubble in $\R^3$ with a cubical inner cell, often created by soap bubble magicians, is stable. 
        }
\end{figure}

\begin{figure}[htbp]
    \centering
     \includegraphics[scale=0.3]{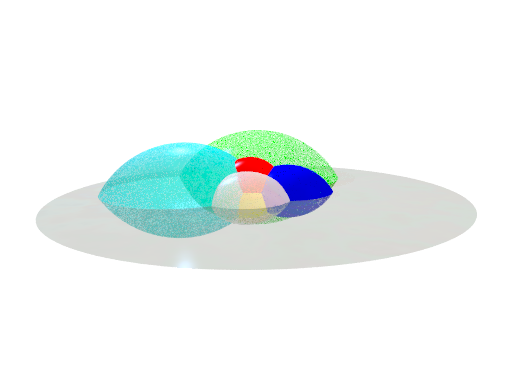}
     \includegraphics[scale=0.3]{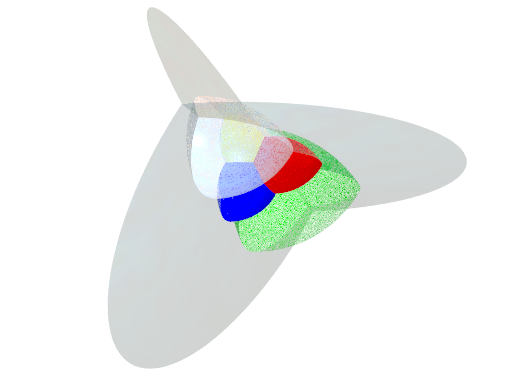}
     \includegraphics[scale=0.25]{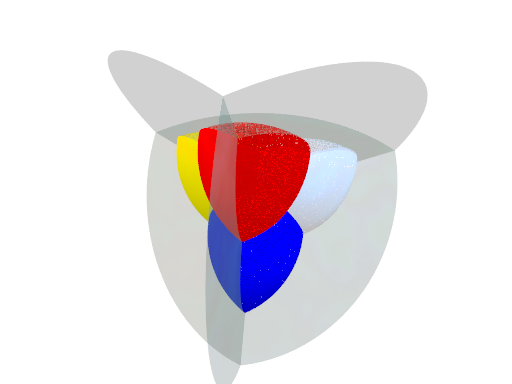}
     \caption{
         \label{fig:nonstandard-partition4} 	Non-standard 8-partitions of $\R^3$ with 2 (left), 3 (middle) and 4 (right) unbounded cells, are all stable. 
        }
\end{figure}

By definition, standard bubbles and partitions satisfy all of the above properties, but there are natural additional examples. 
One just needs to start with a \emph{non-standard} flat regular spherical Voronoi partition in $\S^n$ ($n \geq 3$), and apply a M\"obius transformation to obtain a non-flat partition in $\S^n$, or a stereographic projection onto $\R^n$ or $\HH^n$. For example, an initial non-standard (regular) flat spherical Voronoi $(2n+2)$-partition in $\S^n$ is the following one,
generating the cones over facets of the hypercube in $\R^{n+1}$, namely
\[
 \Omega_i = \set{ p \in \S^n : \argmin_{j=1 \ldots, 2n+2} \scalar{p,\c_j} = \{i\} } ,
\]
for $\c_k = e_k$ and $\c_{k+n+1} = -e_k$, $k=1,\ldots,n+1$. See Figures \ref{fig:nonstandard-partition1} and \ref{fig:nonstandard-partition4} for a depiction of various non-standard bubbles and partitions which satisfy our assumptions and are thus proved to be stable.

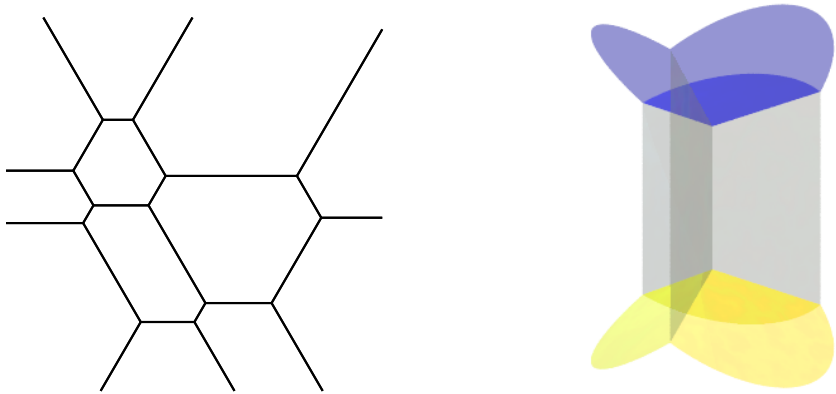
\begin{figure}[htbp]
	\centering
		\begin{tikzpicture}[scale=0.75]				\coordinate (UL) at (-2.8,2.8);
		\coordinate (LL) at (-2.8,-3.8);
		\coordinate (UR) at (3.8,2.8);
		\coordinate (LR) at (3.8,-3.8);
				\coordinate (A) at (0,0);
		\coordinate (M) at (-0.573,0.992);
		\coordinate (N) at (-1.101,0.992);
		\coordinate (P) at (-1.620,0.093);
		\coordinate (I) at (-1.266,-0.520);
		\coordinate (D) at (-0.3,-0.520);
		\draw[thick] (A)--(M)--(N)--(P)--(I)--(D)--(A);
				\coordinate (J) at (-1.447,-0.833);
		\coordinate (K) at (-0.437,-2.582);
		\coordinate (L) at (0.505,-2.582);
		\coordinate (E) at (0.698,-2.247);
		\draw[thick] (I)--(J)--(K)--(L)--(E)--(D)--(I);
				\coordinate (F) at (1.859,-2.247);
		\coordinate (G) at (2.730,-0.739);
		\coordinate (H) at (2.304,0);
		\draw[thick] (E)--(F)--(G)--(H)--(A)--(D)--(E);
				\coordinate (R) at (0.471,2.8);
		\coordinate (Q) at (-2.145,2.8);
		\draw[thick] (R)--(M)--(N)--(Q);
		\coordinate (S) at (-2.8,0.093);
		\draw[thick] (Q)--(N)--(P)--(S);
		\coordinate (U) at (-2.8,-0.833);
		\draw[thick] (S)--(P)--(I)--(J)--(U);
		\coordinate (V) at (-1.140,-3.8);
		\draw[thick] (U)--(J)--(K)--(V);
		\coordinate (W) at (1.208,-3.8);
		\draw[thick] (V)--(K)--(L)--(W);
		\coordinate (Z) at (2.755,-3.8);
		\draw[thick] (W)--(L)--(E)--(F)--(Z);
		\coordinate (A1) at (3.8,-0.739);
		\draw[thick] (Z)--(F)--(G)--(A1);
		\coordinate (T) at (3.8,2.592);
		\draw[thick] (A1)--(G)--(H)--(T);
		\draw[thick] (T)--(H)--(A)--(M)--(R);
	\end{tikzpicture}
	\qquad \hspace{20pt}
     \includegraphics[scale=0.47]{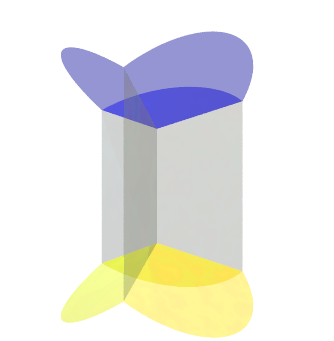}
     \caption{
         \label{fig:Gaussian}
         Two stationary regular flat partitions of $\G^n$: a $12$-partition of $\G^2$ (left) and a $5$-partition of $\G^3$ (right). Both are stable. 
        }
\end{figure}

\medskip
These results provide a partial answer to a question of Kusner \cite[Problem 3]{OpenProblemsInSoapBubbles96}, who asked whether clusters in $\R^3$ with spherical interfaces which meet according to Plateau's laws are necessarily stable. On $\G^n$, we provide a complete answer to the analogous question:
\begin{displayme}
For all $q \geq 2$ and $n \geq 2$, stationary regular flat $q$-partitions in $\G^n$ are stable.
\end{displayme}
A regular partition is called stationary if its interfaces meet in threes at $120^{\circ}$ angles (see below for a precise definition; note that a regular spherical Voronoi partition is automatically stationary). In particular, any collection of line segments (or half-lines) meeting at $120^{\circ}$ angles in the plane delineates a stable partition of $\G^2$ -- see Figure \ref{fig:Gaussian}.

\section{Geometric Measure Theory} \label{sec:GMT}

A powerful starting point for any investigation into multi-bubble isoperimetric problems is provided by geometric measure theory. Using a simple compactness argument for sets of finite perimeter (or functions of bounded variation), it is elementary to show the existence of a minimizing cluster on a space with finite total mass $V_\mu(M) < \infty$. However, this is already non-trivial on $\R^n$, where this was established by Almgren \cite{AlmgrenMemoirs}, who also showed that all finite-volume cells must be bounded; his argument was extended by Morgan \cite{MorganBook5Ed} to any co-compact manifold (where an orbit of some compact set under all isometries covers the entire manifold).

It was also shown by Almgren that the interfaces $\Sigma_{ij} = \partial^* \Omega_i \cap \partial^* \Omega_j$ are $C^\infty$-smooth embedded $(n-1)$-dimensional manifolds with good local separation (see below). Denote 
\[
\Sigma := \cup_{i=1}^q \partial \Omega_i \quad ,\quad \Sigma^1 := \cup_{1 \leq i < j \leq q} \Sigma_{ij}.
\]
In $\R^3$, Plateau's laws were confirmed by J.~Taylor \cite{Taylor-SoapBubbleRegularityInR3}, who showed that $\Sigma$ around $p \in \Sigma \setminus \Sigma^1$ is locally a $C^{1,\alpha}$ diffeomorphic image of a $\Y \times \R$ or $\T$ singularity (in $\R^2$, only a $\Y$ singularity appears). The $\Y$ and $\T$ cones, already alluded to above, are defined as
\begin{align*}
    \Y &:= \{x \in E^{(2)} : \text{ there exist $i \ne j \in \{1,2,3\}$ with $x_i = x_j = \min_{k \in \{1,2,3\}} x_k$}\} , \\
    \T &:= \{x \in E^{(3)} : \text{ there exist $i \ne j \in \{1,2,3,4\}$ with $x_i = x_j = \min_{k \in \{1,2,3,4\}} x_k$}\},
\end{align*}
where $E^{(k)} := \{x \in \R^{k+1} : \sum_{i=1}^{k+1} x_i = 0 \}$ -- see Figure \ref{fig:simplicial-bubble}.
Note that $\Y$ consists of $3$ half-lines meeting at the origin in $120^\circ$ angles, and that $\T$ consists of $6$ two-dimensional sectors meeting in threes at $120^{\circ}$ angles along $4$ half-lines, which in turn all meet at the origin in $\cos^{-1}(-1/3) \simeq 109^{\circ}$ angles.

In $\R^n$ for $n \geq 4$, the analogous verification of Plateau's laws is a more recent result of Colombo--Edelen--Spolaor \cite{CES-RegularityOfMinimalSurfacesNearCones} (which applies in greater generality than just for isoperimetric minimizers, see also \cite{White-AusyAnnouncementOfClusterRegularity,Simon-Codimension2Regularity,NaberValtorta-MinimizingHarmonicMaps}). 
The degree of smoothness locally around a $\Y$ singularity was improved to $C^\infty$ by  Kinderlehrer--Nirenberg--Spruck \cite{KNS} using elliptic regularity for systems of PDEs \cite{ADN2}.
Since these are all local regularity results, they extend to (weighted) Riemannian manifolds. We summarize these and additional properties in the following definition: 
\begin{definition}[Regularity] \label{def:regularity}
A partition $\Omega$ of $(M^n,g,\mu)$ is called regular if it satisfies the following:
\begin{enumerate}
\item $\Omega$ may and will be modified on a $\mu$-null set so that all of its cells are open, and so that for every $i$, $\overline{\partial^* \Omega_i} = \partial \Omega_i$ and $\mu^{n-1} (\partial \Omega_i \setminus \partial^* \Omega_i) =0$; in particular, $\Sigma = \overline{\Sigma^1}$. \item $\Sigma$ is the disjoint union of $\Sigma^1$ and sets $\Sigma^2,\Sigma^3,\Sigma^4$ satisfying (for some fixed $\alpha > 0$):
\begin{enumerate}
\item $\Sigma^1$ is a locally-finite union of embedded $(n-1)$-dimensional $C^{\infty}$ manifolds, and for every $p \in \Sigma^1$, $\Sigma$ around $p$ is locally $C^\infty$ diffeomorphic to $\{0\} \times \R^{n-1}$.
\item $\Sigma^2$ is a locally-finite union of embedded $(n-2)$-dimensional $C^{\infty}$ manifolds, and for every $p \in \Sigma^2$, $\Sigma$ around $p$ is locally $C^\infty$ diffeomorphic to $\Y \times \R^{n-2}$.
\item $\Sigma^3$ is a locally-finite union of embedded $(n-3)$-dimensional $C^{1,\alpha}$ manifolds, and for every $p \in \Sigma^3$, $\Sigma$ around $p$ is locally $C^{1,\alpha}$ diffeomorphic to $\T \times \R^{n-3}$.
\item $\Sigma^4$ is closed and has locally-finite $\H^{n-4}$ measure. 
\end{enumerate}
\item (Density upper bound) For any compact set $K$ in $M$, there exist constants $\Lambda_K, r_K>0$ so that:
\[ \mu^{n-1} (\Sigma \cap B(p, r)) \leq \Lambda_K r^{n-1} , \quad \forall  p \in \Sigma \cap K \quad \forall r \in (0, r_K).
\] \end{enumerate}
\end{definition}

A (locally) minimizing partition is thus always regular (see \cite{EMilmanNeeman-TripleAndQuadruple,EMilmanXu-Stability} and the references therein).
Regularity implies that every point in $\Sigma^2$ (called the \emph{triple-point set}) belongs to the closure of exactly three cells (as well as to the closure of exactly three interfaces). Consequently, $\Sigma^2$ is the disjoint union of $\Sigma_{ijk}$, the subset of $\Sigma^2$ which belongs to the closures of $\Omega_i$, $\Omega_j$ and $\Omega_k$ (or equivalently $\Sigma_{ij}$, $\Sigma_{jk}$ and $\Sigma_{ki}$), and we denote $\partial \Sigma_{ij} := \cup_{k} \Sigma_{ijk}$. 

Let $\Omega$ be a regular $q$-partition. Define the unit normal field $\n_{ij}$ on the interface $\Sigma_{ij}$ pointing from $\Omega_i$ to $\Omega_j$, as well as the corresponding second fundamental form $\II^{ij}$; by regularity, these extend to $\partial \Sigma_{ij}$. The first variation of (weighted) area of $\Sigma_{ij}$ in the normal direction is given by the (weighted) mean curvature $H_{\Sigma_{ij},\mu}$, defined as $\tr(\II^{ij}) +\scalar{\nabla \log \Psi, \n_{ij}}$. 
A simple first variation argument verifies that a (locally) minimizing partition is \emph{stationary}, namely a critical point (with respect to smooth perturbations) of the functional $\F_{\lambda}(\Omega) := A_\mu(\Omega) - \scalar{\lambda, V_\mu(\Omega)}$ for some vector $\lambda \in \R^q$ of Lagrange multipliers (the physical interpretation of $\lambda_i$ is that of the air-pressure inside the cell $\Omega_i$).

\begin{lemma}[Stationarity] Let $\Omega$ be a regular $q$-partition in $(M,g,\mu)$. If $\Omega$ is stationary with Lagrange multipliers $\lambda \in \R^q$, then:
\begin{enumerate}[(1)]
\item \label{it:stationarity-k}
Young-Laplace law: for all $i < j$, $H_{\Sigma_{ij},\mu}$ is constant and equal to $\lambda_i - \lambda_j$.
\item \label{it:stationarity-n} 
Interfaces meet in threes at $120^{\circ}$ angles: $\n_{ij} + \n_{jk} + \n_{ki} = 0$ on $\Sigma_{ijk}$ for all $i < j < k$. 
\end{enumerate}
Conversely, if $\Omega$ is of locally bounded curvature  (i.e.~$\II_{\Sigma^1}$ is bounded on every compact set) and satisfies \ref{it:stationarity-k} and \ref{it:stationarity-n}, then it is stationary with Lagrange multipliers $\lambda \in \R^q$. 
\end{lemma}

Given a vector-field $X \in C_c^\infty(M)$ supported in the interior of a compact $K \subset M$, the flow along $X$ for time $t$ is denoted by $T_t$, the perturbed partition is denoted by $T_t(\Omega) = (T_t(\Omega_i))_{i=1,\ldots,q}$, and its $m$-th variation of volume and total perimeter are denoted by:
\[
\delta^m_X \Vol := \left . \frac{d^m}{(dt)^m} \right|_{t=0} V_\mu(T_t(\Omega);K) \in E^{(q-1)} ~,~ \delta^m_X \Area := \left .  \frac{d^m}{(dt)^m}  \right |_{t=0} A_\mu(T_t(\Omega);K)  \in \R ,
\]
where $V_\mu(\Omega ; K) = (\mu^n( \Omega_i \cap K))_i$ and $A_\mu(\Omega; K) =  \sum_{i<j} \mu^{n-1}(\Sigma_{ij} \cap K)$. 
Stationarity is the first-variation property that $\delta^1_X \F_\lambda(\Omega) = 0$ for all $X \in C_c^\infty(M)$. Additional vital information is contained in the non-negativity of the second-variation of $\F_\lambda(\Omega)$ under a volume constraint. 

\begin{definition}[Stability and Index Form]
A stationary regular partition $\Omega$ is called stable if for every vector-field $X \in C_c^\infty(M)$:
\[
\delta^1_X \Vol = \vec 0 \;\; \Rightarrow \;\; Q(X) := \delta^2_X \Area  - \scalar{\lambda, \delta^2_X \Vol} \geq 0 . 
\]
The quadratic form $Q$ is called the partition's index form. 
\end{definition}

A (locally) minimizing partition $\Omega$ is thus regular, stationary and stable. An additional useful property is the \emph{infiltration property} \cite{Leonardi-Infiltration}, stating the existence of a constant $\eps > 0$ (depending solely on $n$) so that for any $p \in M^n$ and $i=1,\ldots,q$,
\begin{equation} \label{eq:infiltration}
p \in \overline{\Omega_i} \;\; \Rightarrow \;\;  \liminf_{r \rightarrow 0+} \frac{\H^n(\Omega_i \cap B(p,r))}{\H^n(B(p,r))} \geq \eps .
\end{equation}

One final piece of useful information pertains to the symmetry of a minimizing cluster. 

\begin{definition}[$\S^m$-symmetry]
A partition $\Omega$ of $\M^n$, $\M  \in \{ \R , \S, \HH\}$, is said to have $\S^m$-symmetry ($m \in \{ 0, \ldots,n-1\}$), if there exists a totally-geodesic $\M^{n-1-m} \subset \M^n$ so that each cell $\Omega_i$ is invariant under all isometries of $\M^n$ which fix the points of $\M^{n-1-m}$. 
In particular, $\S^0$-symmetry means invariance under reflection about some totally-geodesic hypersurface $\M^{n-1}$. 
\end{definition}

By a simple application of the Borsuk-Ulam (``Ham-Sandwich'') theorem, for any $k$-cluster in $\M^n$ with $k \leq n$ there exists a hyperplane bisecting all of its finite-volume cells. Consequently, by symmetrizing a minimizing cluster about that hyperplane (reflecting the half with the smaller total perimeter), one can always find a minimizing $k$-cluster with $\S^0$-symmetry whenever $k \leq n$. Moreover, when $k \leq n-1$, an argument of White and Hutchings \cite{Hutchings-StructureOfDoubleBubbles} shows that \emph{any} minimizing $k$-cluster will \emph{necessarily} have $\S^{n-k}$-symmetry. On $\M^n \in \{\R^n, \S^n\}$, the spherical Voronoi structure actually yields the same conclusion for all $k \leq n$ \cite{EMilmanNeeman-TripleAndQuadruple}.

\section{Sketch of proofs} \label{sec:proofs}

In this section, we provide a general overview of the strategy for establishing our isoperimetric results described in Section \ref{sec:partial}. We conclude this section by briefly mentioning some ingredients that go into the stability results from Section \ref{sec:stability}.

Let $\Omega$ be a minimizing $k$-cluster in $(M,g,\mu)$; in particular, it is regular, stationary and stable. 
Let $\Psi = \exp(-W)$ denote the density of $\mu$ with respect to $\vol_g$. 
Given a vector-field $X \in C_c^\infty(M)$, define the scalar-field $f$ to be the tuple $(f_{ij})_{i \neq j}$ where $f_{ij} = X^{\n_{ij}}$ on $\Sigma_{ij} \cup \partial \Sigma_{ij}$ is its normal component.
Note that $f_{ji} = -f_{ij}$ is oriented and that $f_{ij} + f_{jk} + f_{ki} = 0$ on $\Sigma_{ijk}$ by stationarity (so called Kirchhoff-Dirichlet boundary conditions). 
Clearly $\delta^1_X V_\mu$ only depends on $f$, namely 
\[
\delta^1_X \Vol = \delta^1_f \Vol := \Big (\sum_{j \neq i} \int_{\Sigma_{ij}} f_{ij} d\mu^{n-1} \Big )_i \in E^{(k)}.
\]

A less trivial fact is that, under appropriate assumptions, the index-form $Q(X)$, which accounts for the second-variation of $\F_{\lambda}(\Omega)$, only depends on $f$ as well. Note that regularity only ensures that $\Sigma$ is $C^{1,\alpha}$ smooth near $\Sigma^3$, and so the curvature $\II = \II_{\Sigma^1}$ could be blowing up near $\Sigma^3$. However, elliptic boundary regularity for systems of PDEs implies that  $\II$ is in $L^2(\Sigma^1 \cap K)$ and in $L^1(\Sigma^2 \cap K)$ for any compact set $K$ disjoint from $\Sigma^4$ (see \cite[Proposition 5.7]{EMilmanNeeman-GaussianMultiBubble}, \cite[Proposition 2.23]{EMilmanNeeman-TripleAndQuadruple}). This is crucial for justifying various integrations by parts on the incomplete manifold-with-boundary $(\Sigma_{ij}, \partial \Sigma_{ij})$, which are required for establishing the following: 

\begin{proposition} 
Let $X \in C_c^\infty(M)$. Assume that either $X$ is supported away from $\Sigma^4$, or that $\Omega$ has locally bounded curvature. Then $Q(X) = Q^0(f)$, where $f = (f_{ij}) = (X^{\n_{ij}})$ and $Q^0(f)$ is given by any of the following two equivalent expressions (and in particular, all terms below are integrable):
\begin{subequations}
 \label{eq:Q0}
\begin{align} 
\label{eq:Q0LJac} Q^0(f)  := & \sum_{i<j} \biggl ( -\int_{\Sigma_{ij}} f_{ij} L_{Jac} f_{ij} d\mu^{n-1} + \int_{\partial \Sigma_{ij}} (\nabla_{\n_{\partial ij}} f_{ij} - \bar \II {\vphantom{\II}^{\partial ij}} f_{ij}) f_{ij} d\mu^{n-2} \biggr) \\
\label{eq:Q0nabla}  = & \sum_{i<j}  \biggl ( \;\; \int_{\Sigma_{ij}} \brac{|\nabla_{\Sigma^1} f_{ij}|^2 - (L_{Jac} 1) f_{ij}^2 } d\mu^{n-1} - \int_{\partial \Sigma_{ij}} \bar \II {\vphantom{\II}^{\partial ij}} f_{ij}^2 d\mu^{n-2} \biggr ) . 
\end{align}
\end{subequations}
\end{proposition}
\noindent
Here $L_{Jac}$ denotes the (weighted) Jacobi operator, defined as
\[ L_{Jac} f_{ij} := \Delta_{\Sigma^1,\mu} f_{ij} + (\Ric_{M,\mu}(\n_{ij},\n_{ij}) + \snorm{\II^{ij}}^2) f_{ij} ,
\] where
\[
\Delta_{\Sigma^1,\mu} f_{ij} := \Delta_{\Sigma^1} f_{ij} - \scalar{\nabla_{\Sigma^1} f_{ij} , \nabla_{\Sigma^1} W} 
\]
is the $\mu$-weighted Laplacian on $(\Sigma^1,\mu^{n-1})$, $\Ric_{M,\mu} := \Ric_M + \nabla^2  W$
is the $\mu$-weighted Ricci curvature on $(M,g,\mu)$, and $\norm{\cdot}$ denotes the Hilbert-Schmidt norm. In addition, 
$\n_{\partial ij}$ denotes the outer co-normal to $\Sigma_{ij}$ on $\partial \Sigma_{ij}$, and $\bar \II^{\partial ij}$ is defined on $\Sigma_{ijk}$ as
\[  \bar \II^{\partial ij} := \frac{\II^{ik}(\n_{\partial ik},\n_{\partial ik}) + \II^{jk}(\n_{\partial jk},\n_{\partial jk})}{\sqrt{3}}.
 \]Consequently, stability amounts to the property that for any scalar-field $f$ as above,
\begin{equation} \label{eq:stability}
\delta^1_f \Vol = \vec 0 \;\; \Rightarrow \;\; Q^0(f) \geq 0 . 
\end{equation}

\subsection{Interfaces have constant curvature}

In the first step, which is the most crucial one, we need to show that all of the interfaces $\Sigma_{ij}$ are locally flat in $\G^n$, i.e.~$\II^{ij} \equiv 0$, or 
locally have constant curvature in $\R^n$ and $\S^n$, i.e.~$\II^{ij}_0 := \II^{ij} - \frac{1}{n-1} \tr(\II^{ij}) \Id \equiv 0$. A classical idea from the single-bubble setting is to test stability on some well-chosen scalar-fields $f$, yielding an integral expression which is on one hand non-positive, and on the other non-negative by stability, implying that the integrand must be equal to zero identically, and hence hopefully that $\II=0$ or $\II_0 = 0$. One challenge in extending this to the multi-bubble setting is a technical one, since the scalar-field $f = (f_{ij})$ needs to be approximated in $H^1(\Sigma^1,\mu)$ by $X^{\n_{ij}}$ around the problematic $\Sigma^3$ (where there isn't enough regularity and curvature may be blowing up) and away from $\Sigma^4$. A more substantial challenge is the boundary term in (\ref{eq:Q0LJac}) or (\ref{eq:Q0nabla}) involving the curvature $\bar \II {\vphantom{\II}^{\partial ij}}$, which may have an \emph{arbitrary sign}. 
A useful idea is therefore to average over a family of scalar-fields $f^\alpha$, selected so that $\delta^1_{f^{\alpha}} \Vol = \vec 0$ for all $\alpha$, and so that in expectation $\E_{\alpha} Q^0_{\mrm{bd}}(f^\alpha) = 0$, where $Q^0_{\mrm{bd}}(f)$ denotes the boundary term in (\ref{eq:Q0LJac}). 
The simplest way is to define $f^a_{ij} = a_{ij} \Phi$ on $\Sigma_{ij}$ for some globally defined function $\Phi \in C^\infty(M)$, where $a_{ij} = a_i - a_j$ and $a \in \R^{k+1}$. Then $\tr( a \mapsto Q^0_{\mrm{bd}}(a_{ij} \Phi)) = 0$ because 
$\n_{\partial ij} + \n_{\partial jk} + \n_{\partial ki} = 0$ and $\bar \II {\vphantom{\II}^{\partial ij}} + \bar \II {\vphantom{\II}^{\partial jk}} + \bar \II {\vphantom{\II}^{\partial ki}}=0$ by stationarity. 
It remains to select a useful $\Phi$, such that $L_{Jac} \Phi$ carries information on $\II$ on one hand, and satisfying \smash{$\int_{\Sigma_{ij}} \Phi d\mu^{n-1} = 0$} for all $i<j$ on the other (to ensure that $\delta^1_{f^a} \Vol = \vec 0$ for all $a$). 

On $\G^n$, when $k \leq n-1$, testing the stability of a cluster under translations easily verifies that necessarily $\Omega$ has a non-trivial product structure $\tilde \Omega \times \R^{n-k}$. A good choice is then to use $\Phi(x) = x_n$, the first eigenfunction of the Gaussian-weighted Laplacian $\Delta_{\R,\gamma^1}$ on the  $\R$ factor in the $n$-th coordinate. If the cluster has no symmetries, we use $f^a_{ij} = a_{ij} - \scalar{\theta(a),\n_{ij}}$, where $\theta(a) \in \R^n$ is a constant translation field chosen to ensure that \smash{$\delta^1_{f^a} \Volg = \vec 0$}, which is always possible when $k=n$. 
In either case, stability verifies that $\II \equiv 0$. 

On $\R^n$ or $\S^n$, we are only able to treat the case when $\Omega$ is $\S^0$-symmetric (invariant under reflection around some $N^{\perp}$), which we may always assume is the case when $k \leq n$, but not when $k=n+1$. A good choice of $\Phi$ is given (again!) by $\Phi(p) = \scalar{p,N}$, but requires additional arguments to conclude that $\II_0 \equiv 0$ on $\R^n$. On $\S^n$ things are more complicated, and we also need to test the stability of $\Omega$ with respect to M\"obius fields, the conformal Killing fields which generate the group of M\"obius automorphism of $\S^n$. Combining the stability from these two families, we deduce that $\II_0 \equiv 0$ (and moreover, that the quasi-centers of all interfaces lie on $N^{\perp}$). 

\subsection{Spherical Voronoi structure}

On $\R^n$ and $\S^n$, there is still quite a bit of work to deduce Theorem \ref{thm:Voronoi} on the spherical Voronoi structure of $\Omega$ and the connectedness of its cells. One can also do this on $\G^n$ (see \cite[Theorem 12.1]{EMilmanNeeman-GaussianMultiBubble}), establishing that $\n_{ij} = \n_i - \n_j$ for each non-empty interface $\Sigma_{ij}$ and some $\{ \n_i \}_{i=1}^{k+1} \subset \R^n$, but this is not as essential in that setting. 

By stereographically projecting from $\R^n$ to $\S^n$, it suffices to treat the latter case. Denote the number of cells of $\Omega$ by $q = k+1$. By passing to the $q'$-partition of $\S^n$ obtained from the connected components (modulo $\S^0$-symmetry) of $\Omega$'s cells, which is still minimizing given the new $q'$ volume constraints, we first assume that all cells are connected (modulo $\S^0$-symmetry) and establish the spherical Voronoi structure in that case; in the next subsection we will use the spherical Voronoi structure to show that all cells must indeed be connected. 

Using an orthogonal projection onto the equatorial plane $\Pi : \S^n \rightarrow N^{\perp}$, we first observe that all projected cells $\Pi \Omega_i$ are convex. Certainly, their boundaries are locally convex around points in $\Pi \Sigma^1$, $\Pi \Sigma^2$ and $\Pi \Sigma^3$ (where locally they are given by the intersection of one, two or three halfplanes), but what about around points in $\Pi \Sigma^4$ (where we have no information)? It turns out that this is irrelevant, thanks to the following extension of a classical local-to-global convexity result of Tietze and Nakajima (corresponding to the case $B = \emptyset$ below) which we establish (see \cite[Proposition 8.7]{EMilmanNeeman-GaussianMultiBubble}):

\begin{proposition}\label{prop:local-to-global-convexity}
    Let $K$ be an open, connected subset of $\R^n$, and let $B$ be a Borel set with $\H^{n-2}(B) = 0$.
    Assume that for every $p \in \partial K \setminus B$ there exists an open neighborhood $U_p$ of $p$ so that
    $K \cap U_p$ is convex. Then $K$ is convex.
\end{proposition}

Once the convexity of all projected cells $\Pi \Omega_i$ is established, their mutual interfaces must also be convex and hence connected. Therefore, each non-empty interface $\Sigma_{ij}$ is a relatively open subset of a single geodesic sphere in $\S^n$ with fixed curvature $\k_{ij}$ and quasi-center $\c_{ij} = \n_{ij} - \k_{ij} p$. Every convex set is the intersection of its supporting halfplanes, and so to establish the spherical Voronoi structure, it remains to show that $\c_{ij} = \c_i - \c_j$ for some collection of $\{\c_i\}_{i=1,\ldots,q} \subset N^{\perp}$ (we already know that $\k_{ij} = \k_i - \k_j = \frac{1}{n-1} (\lambda_i - \lambda_j)$ by the Young-Laplace law). To this end, we consider the two-dimensional simplicial complex whose vertices, edges and triangles are given by $\SS_0 = \{i\}_{i = 1, \ldots,q}$, $\SS_1 = \{ \{i,j\} ; \Sigma_{ij} \neq \emptyset\}$ and $\SS_2 = \{ \{i,j,k\} ; \Sigma_{ijk} \neq \emptyset\}$, and establish that its first cohomology (over any field) vanishes. Since $\c_{ij} + \c_{jk} + \c_{ki} = 0$ for every $\Sigma_{ijk} \neq \emptyset$ by stationarity, it follows that $\c_{ij} = \c_i - \c_j$ for some $\{\c_i\}$ as asserted.

\subsection{Connectedness}

To establish the connectedness of the cells of $\Omega$ when $q \leq n+1$, we first show that every cell which intersects the equator $\S^n \cap N^{\perp}$ (``equatorial cell'') must be connected. Otherwise, we could split a non-connected equatorial cell (say $\Omega_1$) in two, yielding a minimizing $(q+1)$-partition $\Omega'$ in $\S^n$ with $\Omega'_1$ equatorial and $\Omega_{q+1} = \Omega_1 \setminus \Omega'_1$. Consider the (non-smooth yet Lipschitz) scalar-field $f_{ij} = a_{ij} \abs{\scalar{p,N}}$ on the interfaces $\Sigma'_{ij}$ of $\Omega'$. We can ensure by an appropriate choice of $a \in \R^{q+1}$ that $\delta^1_f \Volm(\Omega') = e_1 - e_{q+1}$, and hence \smash{$\delta^1_f \Volm(\Omega) = \vec 0$}. A calculation confirms that \smash{$Q^0_{\Omega}(f) = 0$}, and so it follows by stability (\ref{eq:stability}) and a variational argument that $L_{Jac} f_{ij}$ must be constant on each interface $\Sigma_{ij}$. Elliptic regularity then implies that every $f_{ij}$ must be smooth on $\Sigma_{ij}$, leading to a contradiction, since $a_{1j} \abs{\scalar{p,N}}$ is non-smooth at $\Sigma'_{1j} \cap N^{\perp}$ when $a_{1j} \neq 0$ (and we can ensure that $a_{1j} > 0$ for all $j > 1$ using a strong discrete maximum principle). 

Lastly, we show that when $q \leq n+1$, all cells must be equatorial, and hence connected. The idea is to use the spherical Voronoi structure of the equatorial cells to claim that if there were only $s \leq q-1 \leq n$ equatorial cells, there would not be enough of these to prevent a non-equatorial cell from reaching the equator.

\subsection{Double, triple and quadruple bubbles}

\begin{figure}
    \begin{center}
  \includegraphics[scale=0.09]{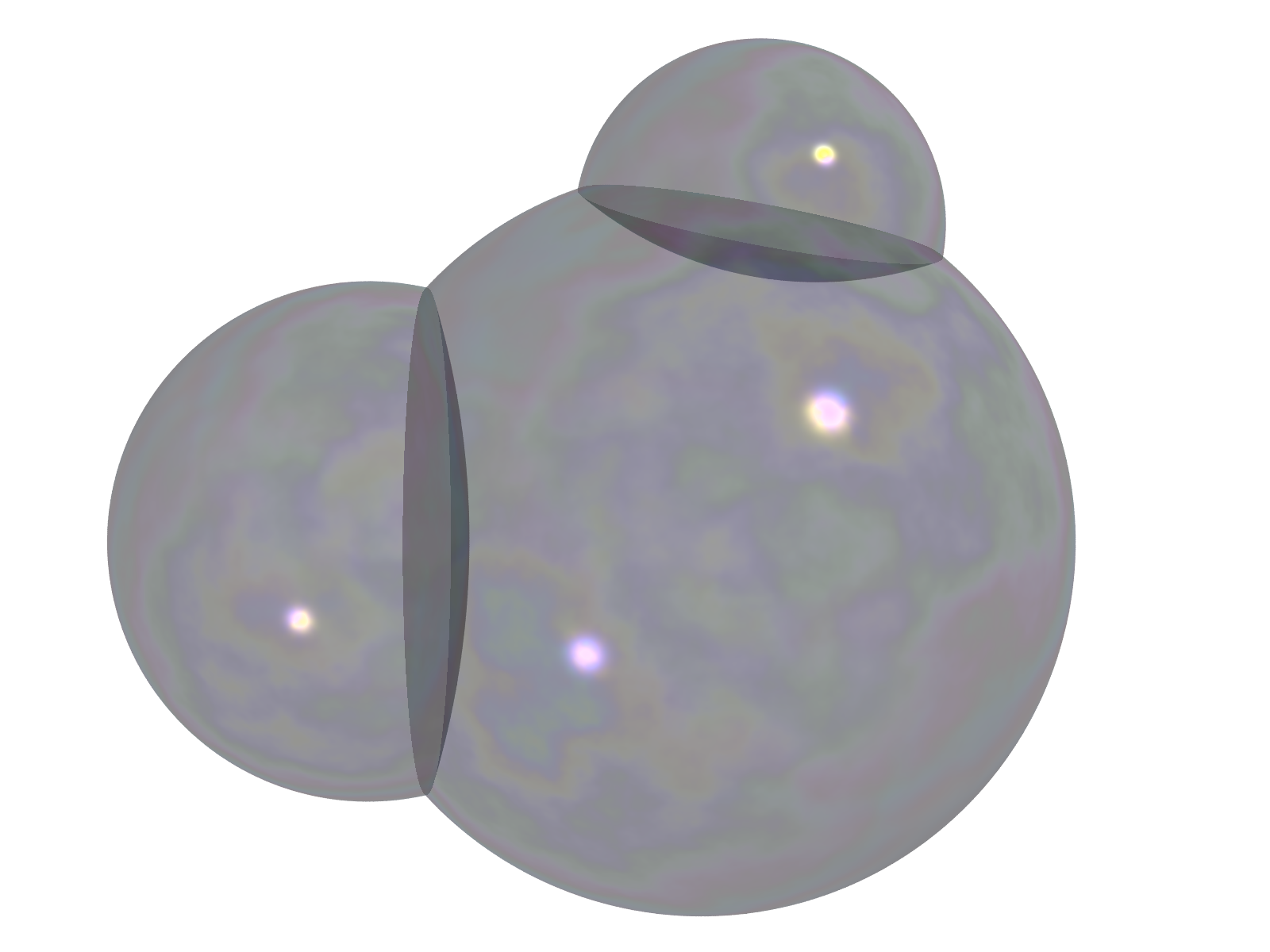}
        \includegraphics[scale=0.09]{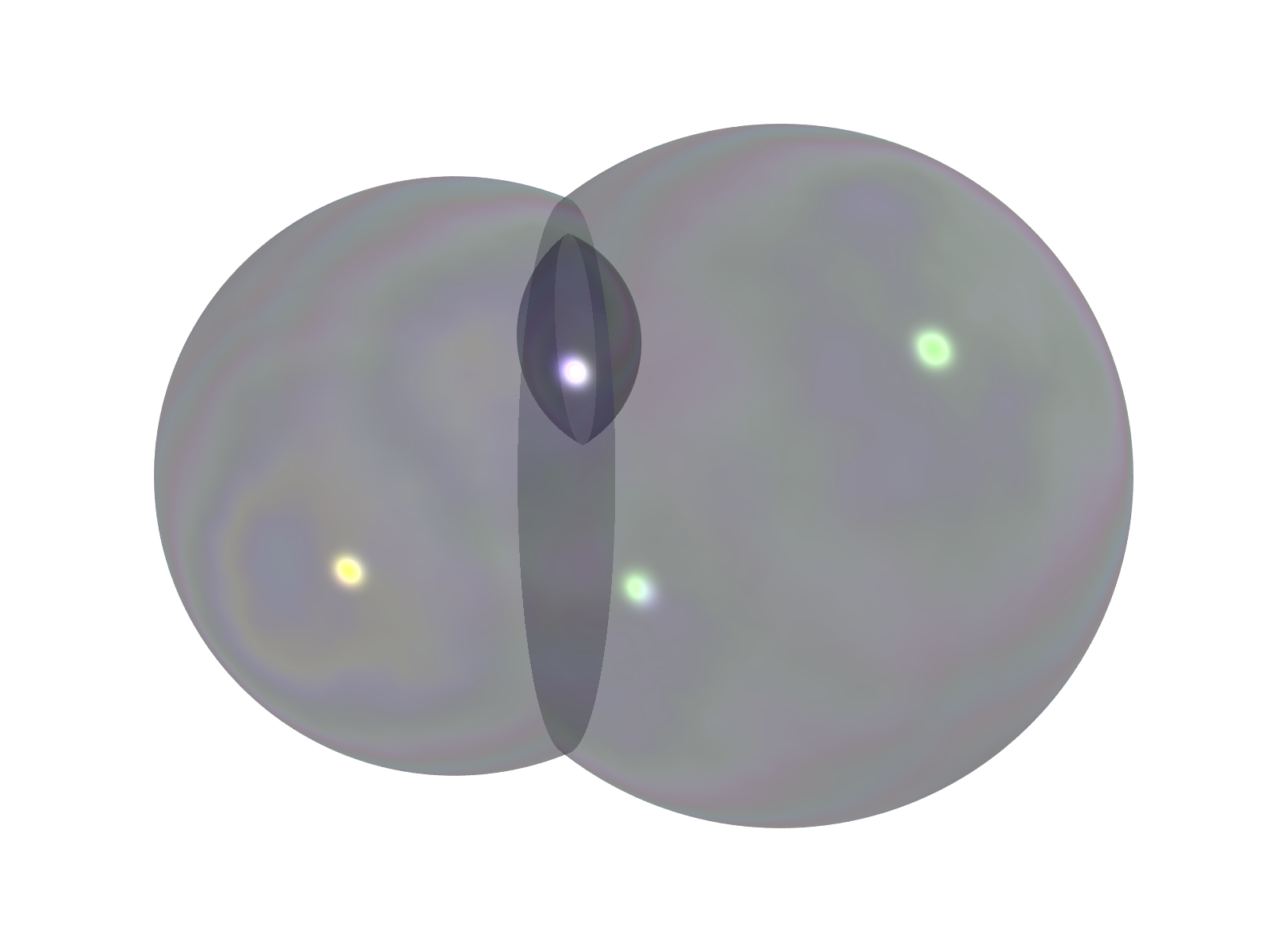}
        \includegraphics[scale=0.08]{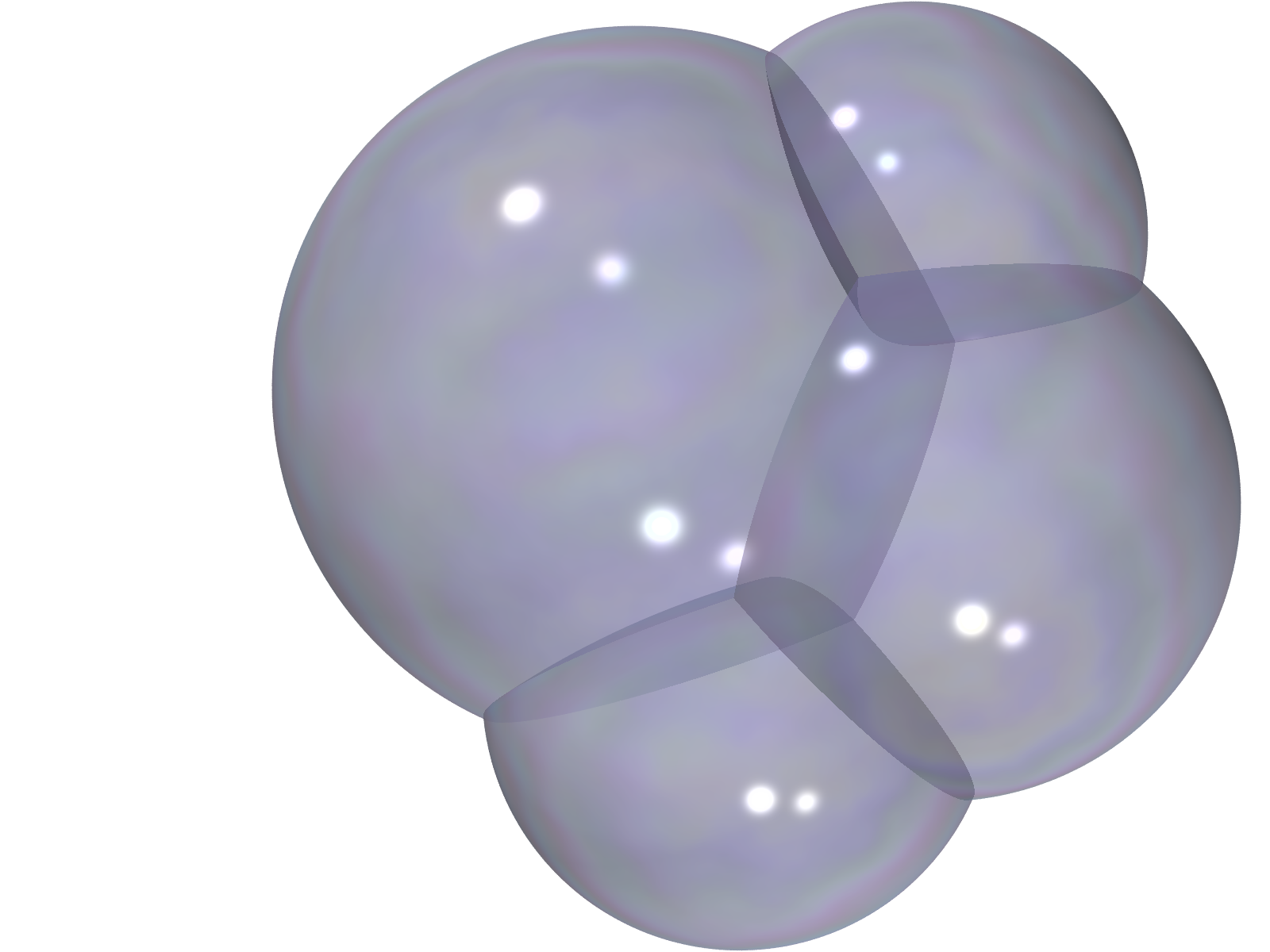}
     \end{center}
     \caption{
         \label{fig:rule-out}
         Non-standard triple-bubbles (left) and quadruple-bubble (right) to be ruled out as global minimizers.  
     }
\end{figure}

So far we've only used the information that $\Omega$ is a \emph{local} minimizer, in the form of stability. To identify the \emph{global} minimizers, we must rule out the other local minimizers, such as disjoint spheres, or the non-standard triple- and quadruple-bubbles depicted in Figure \ref{fig:rule-out}. As already mentioned, when $k \leq n+1$,  standard $k$-bubbles are characterized as those spherical Voronoi $k$-clusters whose interfaces $\Sigma_{ij}$ are non-empty for all $1 \leq i < j \leq k+1$, reducing the problem to showing that the cell-incidence graph of a global minimizer must be the complete graph on $k+1$ vertices.

Assume that a certain bubble configuration is non-rigid, meaning that it is possible to apply some local isometry (like rotation or translation of a certain strict subset of bubbles), preserving all bubble-volumes and total perimeter. A standard idea in geometric measure theory is then to apply the local isometry until a new collision between bubbles occurs, arriving at a contradiction to the known regularity results for isoperimetric minimizing bubbles (given e.g.~in Definition \ref{def:regularity} of regularity or the infiltration property (\ref{eq:infiltration})).

Using this, we are able to show that the incidence graph is not only connected, but in fact 2-connected (removal of any vertex does not disconnect it). This already confirms the double-bubble conjecture on $\R^n$ and $\S^n$, since such a graph on $3$ vertices must be complete. When $3 \leq k \leq 4$, we are also able to show that the minimal degree is at least $3$, confirming the triple-bubble conjectures on $\R^n$ and $\S^n$ ($n \geq 3$). When $k=4$, this leaves only two additional possible graphs on $5$ vertices besides the complete one, one of which is non-rigid and is ruled out as above, and the other corresponds to a configuration arranged in a bubble-ring, which is ruled out by additional geometric arguments. This confirms the quadruple-bubble conjectures on $\R^n$ and $\S^n$ ($n \geq 4$). 

Unfortunately, when the number of bubbles $k$ grows larger, ruling out local minimizers as being global ones using a case-by-case analysis of the cell-incidence graphs on $k+1$ vertices becomes intractable. Fortunately, there is a way to turn local information about minimizers into global information, described next. 

\subsection{Partial differential inequality (PDI) for the isoperimetric profile}

To handle the case of $\G^n$ as well as the quintuple-bubble case on $\S^n$, we examine their corresponding multi-bubble isoperimetric profiles. 
Let us normalize the volume on $\S^n$ so that $V(\S^n) = 1$, and abbreviate $\Delta^{\hspace{-1pt}(k)} = \Delta^{\hspace{-1pt}(k)}[1]$. The tangent spaces of $\Delta^{\hspace{-1pt}(k)}$ are naturally identified with $E^{(k)}$. On $\G^n$, $\mu$ denotes the Gaussian measure and $k \leq n$, whereas on $\S^n$, $\mu$ denotes the Haar probability measure and $k \leq n+1$. 

\begin{definition}[Multi-bubble (model) isoperimetric profile for $\M^n$ ($\M \in \{ \G, \S\}$)]
The multi-bubble isoperimetric profile $\I^{(k)}[\M^n] : \Delta^{\hspace{-1pt}(k)} \to \R_+$ for $k$-clusters in $\M^n$ is defined as:
\[
  \I^{(k)}[\M^n](v) := \inf\{A_\mu(\Omega) : \text{$\Omega$ is a $k$-cluster in $\M^n$ with $V_\mu(\Omega) = v$}\}.
\]
The multi-bubble \emph{model} isoperimetric profile $\I^{(k)}_m[\M^n] : \Delta^{\hspace{-1pt}(k)} \to \R_+$ is defined for $v \in \interior \Delta^{\hspace{-1pt}(k)}$ as:
\[
  \I^{(k)}_m[\M^n](v) := A_\mu(\Omega^m) \text{ where $\Omega^m$ is a standard $k$-bubble in $\M^n$ with $V_\mu(\Omega^m) = v$} .
\]
For $v \in \partial \Delta^{\hspace{-1pt}(k)}$ we define recursively $\I^{(k)}_m[\M^n](v) := \I^{(k-1)}_m[\M^n](v_{-i})$ if $v_i = 0$, where $v_{-i}$ denotes erasing the $i$-th coordinate from $v$ (as there is at least one empty cell). 
\end{definition}
\begin{proposition} 
For $\M \in \{\G, \S\}$, $\I^{(k)}_m[\M^n]$ satisfy the following PDEs on $\Delta^{\hspace{-1pt}(k)}$:
\begin{subequations}
 \label{eq:PDE}
 \begin{align}
\label{eq:PDE-G} & \tr((-\nabla^2 \I^{(k)}_m[\G^n])^{-1})   = 2\I^{(k)}_m[\G^n]  , \\
\label{eq:PDE-S} & \tr\brac{(-\nabla^2 \I^{(k)}_m[\S^n])^{-1} \brac{\frac{1}{2} \Id+ \frac{1}{(n-1)^2} \nabla \I^{(k)}_m[\S^n] \otimes \nabla \I^{(k)}_m[\S^n]}} = \frac{1}{n-1} \I^{(k)}_m[\S^n] .
\end{align}
\end{subequations}
In both cases $-\nabla^2 \I^{(k)}_m[\M^n] > 0$, and hence both model profiles are strictly concave on $\Delta^{\hspace{-1pt}(k)}$, and the above 
 PDEs are elliptic (albeit fully non-linear). \end{proposition}
\begin{remark}
One can similarly define the model isoperimetric profile $\I^{(k)}_m[\R^n] : \Delta^{\hspace{-1pt}(k)}[\infty] \to \R_+$ for $\R^n$ and $k \leq n+1$. However, the homogeneity under scaling of volume and perimeter in $\R^n$ renders the corresponding PDE 
\[
\tr((-\nabla^2 \I^{(k)}_{m}[\R^n])^{-1} \nabla \I^{(k)}_m[\R^n] \otimes \nabla \I^{(k)}_m[\R^n])) = (n-1) \I^{(k)}_m[\R^n] 
\] 
only degenerate elliptic due to the rank one term above, which is why we only consider $\G^n$ and $\S^n$. 
\end{remark}

Our strategy (described below) is to show that the actual profile $\I^{(k)} = \I^{(k)}[\M^n]$ satisfies the corresponding PDE in (\ref{eq:PDE}), in fact as an inequality $\leq$ in the viscosity sense. By induction on $k$, we may assume that $\I^{(k)} = \I^{(k)}_m$ on $\partial \Delta^{\hspace{-1pt}(k)}$, and so by an application of the maximum principle (as our PDEs are elliptic), we deduce that $\I^{(k)} \geq \I^{(k)}_m$. In this manner, global information is propagated from $\partial \Delta^{\hspace{-1pt}(k)}$ to $\interior \Delta^{\hspace{-1pt}(k)}$, verifying that $A_{\mu}(\Omega) \geq \I^{(k)}_m(V_\mu(\Omega))$ for all local minimizers $\Omega$. Since trivially $\I^{(k)} \leq \I^{(k)}_m$, we conclude that $\I^{(k)} = \I^{(k)}_m$. Analysis of the equality case requires some additional work, but in essence, this confirms the $k$-bubble conjecture on $\G^n$ for the entire applicable range $2 \leq k \leq n$. On $\S^n$, we are only able to show that $\I^{(k)}$ satisfies the above PDI for $k \leq 5$, thus avoiding the analysis of the cell-incidence graphs when $k \leq 4$, and resolving the quintuple-bubble case $k=5$ on $\S^n$ ($n \geq 5$). By rescaling and approximately embedding a small cluster in $\R^n$ into $\S^n$, we deduce that a standard quintuple-bubble in $\R^n$ ($n \geq 5$) is an isoperimetric minimizer, resolving the quintuple-bubble case in that setting as well; however, we cannot verify that standard quintuple-bubbles are the only minimizers, as uniqueness is lost in the approximation procedure.

Our goal is thus to obtain sharp upper bounds on $\nabla^2 \I^{(k)}$ (in the viscosity sense); to simplify this sketch, let us assume that $\I = \I^{(k)}$ is smooth. Fix an isoperimetric minimizing $k$-cluster $\Omega$ in $\M^n$, $\M \in \{\S , \G\}$ and $k \leq n$.  On $\S^n$ we already know that $\Omega$ is spherical Voronoi, on $\G^n$ it suffices to know that its interfaces are flat. 
If $T_t$ is the flow on $\M^n$ generated by a vector-field $X$, by definition $\I(V_\mu(T_t(\Omega))) \leq A_\mu(T_t(\Omega))$ with equality at $t=0$. Differentiating twice at $t=0$, we obtain at $v_0 = V_\mu(\Omega) \in \interior \Delta^{\hspace{-1pt}(k)}$:
\[
 \scalar{\nabla \I(v_0),\delta^1_X \Vol} = \delta^1_X \Area  ~,~ (\delta^1_X \Vol)^T \nabla^2 \I(v_0) \delta^1_X \Vol  + \scalar{\nabla \I(v_0),\delta_X^2 \Vol} \leq \delta_X^2 \Area .
\]
Since $\delta^1_X \Area = \scalar{\lambda,\delta^1_X \Vol}$ by stationarity (and since $\delta^1_X \Vol \in E^{(k)}$ may be chosen arbitrarily), we see that $\lambda = \nabla \I(v_0)$, and after rearrangement, we deduce:
\begin{equation} \label{eq:Qinq}
(\delta^1_X \Vol)^T \nabla^2 \I(v_0) \delta^1_X \Vol \leq Q(X) . 
\end{equation}
Note that this recovers stability: $\delta^1_X \Vol = \vec 0 \; \Rightarrow \; 0 \leq Q(X)$. The challenge now is to find a rich-enough family of vector-fields $X$ which will yield a sharp upper bound on $\nabla^2 \I(v_0)$ via (\ref{eq:Qinq}). 

\subsection{Which fields to use?}

The first natural idea is to use the family of vector-fields $X$ which generate the conjectured minimizers -- translation vector-fields on $\G^n$ and M\"obius vector-fields on $\S^n$. This works perfectly well when the minimizing $k$-cluster $\Omega$ is ``full-dimensional'' (in ``general position''), namely when the normals $\{ \n_{ij} = \n_i - \n_j \}$ to the interfaces $\{ \Sigma_{ij} \}$ in $\G^n$ or the corresponding quasi-centers $\{\c_{ij} = \c_i - \c_j\}$ in $\S^n$ span a $k$-dimensional linear subspace $\mathcal{N}$. However, we do not know how to \emph{a priori} exclude the possibility that the cells of $\Omega$ have arranged in some lower-dimensional configuration,  i.e.~that $\dim \mathcal{N} < k$,  in which case these generating fields will only give us access to a lower-dimensional minor of the $k$-dimensional quadratic form $\nabla^2 \I(v_0)$. This absurd possibility is the crux of the difficulty in extending our results on $\S^n$ (and by approximation, $\R^n$) to the entire range $k \leq n$. Consequently, we need to test a more general family of vector-fields in (\ref{eq:Qinq}). 

Of course, the tightest inequality in (\ref{eq:Qinq}) is obtained when $X$ is a minimizer of $Q(X)$ under the volume constraint $\delta^1_X \Vol = \delta v$ for a given $\delta v \in E^{(k)}$.
Expressing everything in terms of the normal components $f = (f_{ij})$, we are looking for minimizers of $Q^0(f)$ on the affine subspace $\smash{\{\delta^1_f \Vol = \delta v\}}$  in an appropriate Sobolev space $H^1(\Sigma^1,\mu^{n-1})$. On $\S^n$, we do not know how to guarantee the existence of minimizers for every $\delta v \in E^{(k)}$, since we do not know in general whether $Q^0(f)$ is necessarily bounded below on every affine subspace $\smash{\{\delta^1_f \Vol = \delta v\}}$. However, if a minimizer $f$ does exist, a simple variational argument verifies that it would necessarily be a \emph{conformal Jacobi field}, namely satisfy $L_{Jac} f_{ij} = a_{ij} := a_i - a_j$ on every $\Sigma_{ij}$ for some $a \in E^{(k)}$, and belong to an appropriate domain $\D_{con}$ of Sobolev functions satisfying \emph{conformal boundary conditions} on $\Sigma^2$. These boundary conditions involve a mixture of Kirchhoff-Dirichlet and Robin conditions where $3$ interfaces meet, ensuring that $Q^0_{\mrm{bd}}(f) = 0$. On the domain $\D_{con}$ the Jacobi operator $L_{Jac}$ is self-adjoint and Fredholm in $L^2(\Sigma^1,\mu^{n-1})$, and its spectral theory, which is reminiscent of the one for quantum graphs, plays an important role below.

On $\G^n$, we have an explicit description of a $k$-dimensional family of conformal Jacobi fields --- these are given by piecewise constant scalar-fields $f^a = (a_{ij}) \in \D_{con}$ for $a \in E^{(k)}$, since $L_{Jac} 1 = 1$ (as $\II = 0$ and $\Ric_{\R^n,\gamma^n} \equiv 1$) and hence $L_{Jac} a_{ij} = a_{ij}$ on $\Sigma_{ij}$. Note that the normal components $(\scalar{\theta,\n_{ij}})$ of a constant translation field $X \equiv \theta$ are of this form, but if the cluster is lower-dimensional, not \emph{all} piecewise constant scalar-fields $(a_{ij})$ would be of the form $(\scalar{\theta, \n_{ij}})$ (due to linear dependencies). Denoting
\[
L_\gamma := \sum_{i<j} \gamma^{n-1}(\Sigma_{ij}) (e_i - e_j) \otimes (e_i - e_j) ,
\vspace{-5pt}
\]
it is easy to see that $L_{\gamma}$ is positive-definite on $E^{(k)}$ whenever $V_{\gamma^n}(\Omega) \in \interior \Delta^{\hspace{-1pt}(k)}$. 
Computing, we have
\[
\delta^1_{f^a} V_{\gamma^n} = \Big ( \sum_{j \neq i} a_{ij} \gamma^{n-1}(\Sigma_{ij}) \Big  )_i = L_\gamma a \quad,\quad Q^0(f^a) =  - \scalar{L_{Jac} f^a , f^a}_{\Sigma^1,\gamma^{n-1}} = - a^T L_\gamma a .
\vspace{-5pt}
\]
By constructing smooth vector-fields $X$ so that $(X^{\n_{ij}})$ approximates $(a_{ij})$ in $H^1(\Sigma^1,\gamma^{n-1})$ (and in some additional metrics), and plugging the above computations into (\ref{eq:Qinq}) for all $a \in E^{(k)}$, we thus confirm the desired PDI (\ref{eq:PDE-G}) for $\I = \I^{(k)}[\G^n]$:
\[
L_\gamma \nabla^2 \I L_\gamma \leq - L_\gamma \; \Rightarrow \; 0 < (-\nabla^2 \I)^{-1} \leq L_\gamma \; \Rightarrow \; \tr((-\nabla^2 \I)^{-1})   \leq \tr(L_\gamma) = 2 \sum_{i<j} \gamma^{n-1}(\Sigma_{ij}) = 2 \I . 
\]

On $\S^n$, the normal component of a M\"obius field is a conformal Jacobi field, but again, if the cluster $\Omega$ is lower-dimensional, these will not span a $k$-dimensional family. In contrast to the case of $\G^n$, we generally have no explicit description of all conformal Jacobi fields for $\Omega$. However, using the Fredholm alternative, we are able to show that a solution $f^a \in \D_{con}$ to $L_{Jac} f^a_{ij} = (n-1) a_{ij}$ always exists for all $a \in E^{(k)}$. Using this, we can already determine the index of the quadratic form $Q^0$ and thereby establish the concavity of the isoperimetric profile \cite{EMilmanNeeman-QuintupleBubble}, answering a question of Heppes \cite[Problem 4]{OpenProblemsInSoapBubbles96}:
\begin{theorem} 
Let $k \leq n$. Then for any minimizing $k$-cluster $\Omega$ in $\S^n$ with $V(\Omega) \in \interior \Delta^{\hspace{-1pt}(k)}$, $(-L_{Jac},\D_{con})$ has exactly $k$ negative eigenvalues. As a consequence, $\I^{(k)}[\S^n]$ is strictly concave on $\Delta^{\hspace{-1pt}(k)}$, and $\I^{(k)}[\R^n]$ is concave on $\Delta^{\hspace{-1pt}(k)}[\infty]$. 
\end{theorem}

\subsection{Trace identity} \label{subsec:trace}
 
 Given a minimizing $k$-cluster in $\S^n$ ($k \leq n$), let $\FF$ denote the linear operator mapping $a \in E^{(k)}$ to $\delta^1_{f^a} \Volm \in E^{(k)}$, 
  where recall $f^a \in \D_{con}$ is the conformal Jacobi field solving $L_{Jac} f^a_{ij} = (n-1) a_{ij}$ on every $\Sigma_{ij}$; it turns out that $\FF$ is well-defined, symmetric, positive semi-definite, and that $Q^0(f^a) = - (n-1) a^T \FF a$. Testing these conformal Jacobi fields in (\ref{eq:Qinq}) as in the Gaussian setting, with $\FF$ now playing the role of $L_\gamma$, the task of verifying the sharp PDI (\ref{eq:PDE-S}) for $\I^{(k)}[\S^n]$ reduces to establishing the following trace-identity for $\FF$ (even just as an inequality $\leq $):
\begin{equation} \label{eq:trace-identity}
\tr\brac{ \FF \brac{\frac{1}{2} \Id + \k \otimes \k}} = \H^{n-1}(\Sigma^1) ,
\end{equation}
where $\k = (\k_1,\ldots,\k_{k+1}) = \frac{\lambda}{n-1} \in E^{(k)}$ are the curvature parameters of $\Omega$ from its spherical Voronoi description.

We can show that (\ref{eq:trace-identity}) holds in a variety of scenarios, certainly if $\Omega$ is full-dimensional or M\"obius-flat, but we were unable to verify (\ref{eq:trace-identity}) in general. We are able to verify (\ref{eq:trace-identity}) for a certain relaxation $\FF_0$ of $\FF$, constructed as the limit $\lim_{t \rightarrow 0} \FF_t$ where $\FF_t$ are the corresponding operators associated to a one-parameter family of conformally perturbed clusters $\Omega_t$. However, we could not verify (in general) that $\FF = \FF_0$, i.e.~that $\FF(\lim_{t \rightarrow 0} \Omega_t) = \lim_{t \rightarrow 0} \FF(\Omega_t)$, which would follow if $\FF(\Omega)$ were continuous in $\Omega$ (!). A confirmation of this innocent-looking PDE question would immediately allow us to extend our results on $\S^n$ and $\R^n$ from the quintuple case to general $k \leq n$. 
 
\subsection{Quintuple bubble}
 
When $k \leq 5$ on $\S^n$, we are able to show that either the minimizing $k$-cluster $\Omega$ satisfies a certain generalized M\"obius-flatness condition, in which case $\FF = \FF_0$ and (\ref{eq:trace-identity}) holds, or else it satisfies a higher-order extension of Plateau's laws, meaning that all tangent cones to $\Sigma$ are cones over (regular) simplices. In the latter case, we are able to slightly perturb $\Omega$ into a full-dimensional cluster without creating any new interfaces. Contrary to other deformations described above, which were local isometries requiring some non-rigidity, this deformation is neither, and yet is guaranteed not to alter the cells' volumes nor total perimeter. We thus obtain a different minimizing cluster for the same volume constraints, which is now full-dimensional, and therefore satisfies (\ref{eq:trace-identity}). This establishes the sharp PDI (\ref{eq:PDE-S}) for $\I^{(k)}[\S^n]$ when $k \leq \min(5,n)$, resolving in particular the quintuple-bubble case when $n \geq 5$. 

\subsection{Establishing stability}

We conclude this section with a few comments on how to derive the results described in Section \ref{sec:stability}. 
Note that stability (\ref{eq:stability}) is a Poincar\'e-type inequality on the collection of interfaces $\Sigma^1 = \cup_{i<j} \Sigma_{ij}$ and their boundaries $\Sigma^2 = \cup_{i<j} \partial \Sigma_{ij}$. Consequently, to establish the stability of a partition, we extend to the multi-bubble setting some known tools from the single weighted manifold-with-boundary setting, namely the $L^2$ Bochner method and Brascamp-Lieb--type inequalities. These methods typically require some control over the curvature of the boundaries $\partial \Sigma_{ij}$, which we do not have unless all of the interfaces are flat (like in $\G^n$). To address this, we are able to construct, at least for regular M\"obius-flat spherical Voronoi partitions of a model space $\M^n \in \{ \R^n, \S^n, \HH^n\}$, a potential function which conformally flattens all boundaries $\partial \Sigma_{ij}$. 
We then incorporate this conformally flattening potential into a conjugated multi-bubble Brascamp-Lieb inequality, which extends the single-bubble one due to Huang--Zhu \cite{HuangMaZhu-ReillyFormula,HuangZhu-ConjugatedBrascampLieb}.

\section{Open problems}

We conclude this survey with a list of open problems (extending the classical list \cite{OpenProblemsInSoapBubbles96} of open problems in soap bubble geometry):
\begin{enumerate}
\item Establish the case $k=n+1$ ($n \geq 3$) in Sullivan's conjectures, and in particular the quadruple-bubble case on $\R^3$. Our results are currently restricted to $k \leq n$, because we need an initial $\S^0$-symmetry to find a rich-enough family of volume-preserving perturbations with which to test stability, but 
we believe that Theorem \ref{thm:Voronoi} should also hold for $k = n+1$. If this is too challenging, show that a standard $(n+2)$-partition in $\R^n$ with at least 2 infinite cells is locally minimizing in the sense of \cite{ABV-LensClusters} -- this might be easier, since the (at most $n$) finite-volume cells will have an $\S^0$-symmetry. Note that in the Gaussian setting $\G^n$, we are actually able to handle the maximal case $k=n$ without any symmetry, by appropriately offsetting our perturbations.
\item Extend the results described above to the hyperbolic setting on $\HH^n$. The double-bubble conjecture has been established on $\HH^2$ \cite{CottonFreeman-DoubleBubbleInSandH}, but remains wide open on $\HH^n$ for $n \geq 3$. It should be possible to establish the triple-bubble conjecture on $\HH^2$ by mimicking the arguments of \cite{Wichiramala-TripleBubbleInR2} and \cite{ Lawlor-TripleBubbleInR2AndS2} in $\R^2$ and $\S^2$, respectively. 
However, the higher-dimensional case seems out of reach of the approach in \cite{EMilmanNeeman-TripleAndQuadruple} without some new ingredient, since one of the inequalities deduced from stability on $\S^n$ goes in the wrong direction on $\HH^n$. The main challenge is to establish that the interfaces of a minimizing cluster in $\HH^n$ are generalized spheres.
 \item Kusner's question \cite[Problem 3]{OpenProblemsInSoapBubbles96} asks whether clusters in $\R^3$ with spherical interfaces which meet according to Plateau's laws are necessarily stable; a natural extension to $\R^n$ would be to require that all tangent cones are area-minimizing, or just that the partition is regular. A related conjecture of Morgan \cite[Conjecture 2.4]{Morgan-StrictCalibrations} asserts that any stationary regular partition in $\R^n$, all of whose tangent cones are strictly area-minimizing, is locally area-minimizing under a volume constraint in a small-enough ball (see \cite[Theorem 2.1]{Morgan-StrictCalibrations} for a proof of this for $C^1$ deformations in a small ball where 3 cells meet, and \cite[Theorem 5.1]{LawlorMorgan-TripleJunctions} for arbitrary deformations of 3 minimal surfaces in a small ball).
\item Obtain a characterization of stable (stationary regular) $q$-partitions $\Omega$ of $\G^n$ for all $q \geq 3$ 
(for a characterization in the case $q=2$, see \cite{McGonagleRoss:15} on $\G^n$, and \cite{Alexandrov-MethodOfMovingPlanes,BarbosaDoCarmo-StabilityInRn,BarbosaDoCarmoEschenburg-StabilityInManifolds} on $\R^n$, $\S^n$ and $\HH^n$). 
Our results from \cite{EMilmanNeeman-GaussianMultiBubble} imply that whenever a stable $\Omega$ has a product structure $\Omega' \times \R$, and in particular, whenever $q \leq n$, then necessarily $\Omega$ has flat interfaces; the same also holds for $q = n+1$ using a separate argument which does not require symmetry.
Conversely, our results from \cite{EMilmanXu-Stability} imply (modulo technicalities) that a stationary regular $q$-partition with flat interfaces is necessarily stable  for any $q \geq 2$. Consequently, the remaining question is whether a stable $q$-partition with $q \geq n+2$ and with no product structure necessarily has flat interfaces -- this is open even on $\G^2$. Note that the analogous question in $\R^n$ is apparently false, since even on $\R^3$, computer simulations suggest that a minimizing $6$-bubble may not have spherical interfaces \cite{SullivanOldSurvey}.

\item Given a minimizing (spherical Voronoi) $k$-cluster in $\S^n$, $k \leq n$, show that the positive semi-definite operator $\FF$ from Subsection \ref{subsec:trace} is in fact always positive-definite. Equivalently, show that $Q^0(f)$ is always bounded below on the affine subspace $\{ f \in \D_{con} : \delta^1_f \Volm = \delta v\}$ for all $\delta v \in E^{(k)}$. 
\item Given a minimizing (spherical Voronoi) $k$-cluster in $\S^n$, $k \leq n$, symmetric under reflection about $N^{\perp}$, establish the trace identity (\ref{eq:trace-identity}). This would extend the confirmation of the multi-bubble isoperimetric conjecture  on $\S^n$ (and without uniqueness, on $\R^n$) from $k \leq 5$ to any $k \leq n$. Several concrete ways for establishing  (\ref{eq:trace-identity}) are as follows (see \cite{EMilmanNeeman-QuintupleBubble}):
\begin{enumerate}
\item Show that $\FF = \FF_0$, where $\FF_0 = n \sum_{i<j} \int_{\Sigma_{ij}} \scalar{p,N}^2 d\H^{n-1}(p) (e_i - e_j) \otimes (e_i - e_j)$. 
\item Establish the continuity $\FF(\lim_{t \rightarrow 0} \Omega_t) = \lim_{t \rightarrow 0} \FF(\Omega_t)$, where $\Omega_t = T_t(\Omega)$ and $\{T_t\}$ denote the one-parameter M\"obius automorphisms of $\S^n$ generated by the (conformal Killing) M\"obius vector-field $N - \scalar{N,p} p$. 
\item Show that $\Sigma_{ij} = \emptyset \;\; \Rightarrow \;\; \FF_{ij} = 0$ for all $1 \leq i < j \leq k+1$. 
\end{enumerate}
\item In order to remove the technical assumptions we require from our test functions in \cite{EMilmanXu-Stability} when testing stability, one would need to extend the Agmon--Douglis--Nirenberg theory of boundary regularity for systems of elliptic PDEs \cite{ADN2} from half-planes $\{ x \in \R^n : x_1 \geq 0 \}$ to convex sectors $\{ x \in \R^n : x_1 \geq 0 , x_1 + a x_2 \geq 0 \}$ ($a \neq 0$). Indeed, the six interfaces $\{\Sigma_{ij}\}_{\{i,j\} \subset \{a,b,c,d\}}$ meeting at a $\T$-type singularity on $\Sigma^3$ locally look like sectors with aperture angle  $\cos^{-1}(-1/3) \simeq 109^{\circ}$. 
 \end{enumerate}

\section*{Acknowledgments.}
The research leading to these results is part of a project that has received funding from the European Research Council (ERC) under the European Union's Horizon 2020 research and innovation programme (grant agreement No 101001677). It is a pleasure to warmly thank Francesco Maggi for his comments on a preliminary version.

\bibliographystyle{plain}
\bibliography{ICM-2026-Bib}

@preamble{ "\def\cprime{$'$} \def\textasciitilde{$\sim$}"
}

@article{Douglas1931-PlateauProblem,
  author =        {Douglas, J.},
  journal =       {Trans. Amer. Math. Soc.},
  number =        {1},
  pages =         {263--321},
  title =         {Solution of the problem of {P}lateau},
  volume =        {33},
  year =          {1931},
  doi =           {10.2307/1989472},
  issn =          {0002-9947},
  url =           {https://doi-org.ezlibrary.technion.ac.il/10.2307/1989472},
}

@article{Rado1930-PlateauProblem,
  author =        {Rad\'{o}, T.},
  journal =       {Math. Z.},
  number =        {1},
  pages =         {763--796},
  title =         {The problem of the least area and the problem of
                   {P}lateau},
  volume =        {32},
  year =          {1930},
  doi =           {10.1007/BF01194665},
  issn =          {0025-5874},
  url =           {https://doi-org.ezlibrary.technion.ac.il/10.1007/
                  BF01194665},
}

@book{MorganBook5Ed,
  author =        {Morgan, F.},
  edition =       {Fifth},
  note =          {A beginner's guide, Illustrated by James F. Bredt},
  pages =         {viii+263},
  publisher =     {Elsevier/Academic Press, Amsterdam},
  title =         {Geometric measure theory},
  year =          {2016},
  isbn =          {978-0-12-804489-6},
}

@incollection{David-ShouldWeSolvePlateau,
  author =        {David, G.},
  booktitle =     {Advances in analysis: the legacy of {E}lias {M}.
                   {S}tein},
  pages =         {108--145},
  publisher =     {Princeton Univ. Press, Princeton, NJ},
  series =        {Princeton Math. Ser.},
  title =         {Should we solve {P}lateau's problem again?},
  volume =        {50},
  year =          {2014},
}

@article{KMS-PlateauAsLimitOfCapillarity,
  author =        {King, D. and Maggi, F. and Stuvard, S.},
  journal =       {Comm. Pure Appl. Math.},
  number =        {5},
  pages =         {895--969},
  title =         {Plateau's problem as a singular limit of capillarity
                   problems (revised)},
  volume =        {75},
  year =          {2022},
  doi =           {10.1002/cpa.22048},
  issn =          {0010-3640},
  url =           {https://doi-org.ezlibrary.technion.ac.il/10.1002/cpa.22048},
}

@unpublished{MaggiNovackRestrepo-PlateauBorders,
  author =        {Maggi, F. and Novack, M. and Restrepo, D.},
  note =          {arXiv:2310.20169},
  title =         {Plateau borders in soap films and {G}auss'
                   capillarity theory},
  year =          {2023},
}

@book{MaggiBook,
  author =        {Maggi, F.},
  pages =         {xx+454},
  publisher =     {Cambridge University Press, Cambridge},
  series =        {Cambridge Studies in Advanced Mathematics},
  title =         {Sets of finite perimeter and geometric variational
                   problems: an introduction to {G}eometric {M}easure
                   {T}heory},
  volume =        {135},
  year =          {2012},
  isbn =          {978-1-107-02103-7},
  url =           {https://doi.org/10.1017/CBO9781139108133},
}

@article{Steiner-1842,
  author =        {Steiner, J.},
  journal =       {J. Reine Angew. Math.},
  pages =         {93--152},
  title =         {Sur le maximum et le minimum des figures dans le
                   plan, sur la sph\`ere et dans l'espace en
                   g\'{e}n\'{e}ral. {P}remier m\'{e}moire},
  volume =        {24},
  year =          {1842},
  doi =           {10.1515/crll.1842.24.93},
  issn =          {0075-4102},
  url =           {https://doi-org.ezlibrary.technion.ac.il/10.1515/
                  crll.1842.24.93},
}

@book{Schwarz-1890,
  author =        {Schwarz, H. A.},
  note =          {Nachdruck in einem Band der Auflage von 1890},
  pages =         {},
  publisher =     {Chelsea Publishing Co., Bronx, NY},
  title =         {Gesammelte mathematische {A}bhandlungen. {B}and {I},
                   {II}},
  year =          {1972},
}

@article{SchmidtIsopOnModelSpaces,
  author =        {Schmidt, E.},
  journal =       {Math. Nachr.},
  pages =         {81--157},
  title =         {Die {B}runn-{M}inkowskische {U}ngleichung und ihr
                   {S}piegelbild sowie die isoperimetrische
                   {E}igenschaft der {K}ugel in der euklidischen und
                   nichteuklidischen {G}eometrie. {I}},
  volume =        {1},
  year =          {1948},
}

@book{BuragoZalgallerBook,
  address =       {Berlin},
  author =        {Burago, Yu. D. and Zalgaller, V. A.},
  publisher =     {Springer-Verlag},
  series =        {Grundlehren der Mathematischen Wissenschaften
                   [Fundamental Principles of Mathematical Sciences]},
  title =         {Geometric inequalities},
  volume =        {285},
  year =          {1988},
}

@article{SudakovTsirelson,
  author =        {Sudakov, V. N. and
                   Cirel{\cprime}son [Tsirelson], B. S.},
  journal =       {Zap. Nau\v cn. Sem. Leningrad. Otdel. Mat. Inst.
                   Steklov. (LOMI)},
  note =          {Problems in the theory of probability distributions,
                   II},
  pages =         {14--24, 165},
  title =         {Extremal properties of half-spaces for spherically
                   invariant measures},
  volume =        {41},
  year =          {1974},
}

@article{Borell-GaussianIsoperimetry,
  author =        {Borell, Ch.},
  journal =       {Invent. Math.},
  pages =         {207-216},
  title =         {The {B}runn--{M}inkowski inequality in {G}auss
                   spaces},
  volume =        {30},
  year =          {1975},
}

@article{CarlenKerceEqualityInGaussianIsop,
  author =        {Carlen, E. A. and Kerce, C.},
  journal =       {Calc. Var. Partial Differential Equations},
  number =        {1},
  pages =         {1--18},
  title =         {On the cases of equality in {B}obkov's inequality and
                   {G}aussian rearrangement},
  volume =        {13},
  year =          {2001},
}

@incollection{SullivanOldSurvey,
  author =        {Sullivan, J. M.},
  booktitle =     {Foams and emulsions ({C}arg\`ese, 1997)},
  pages =         {379--402},
  publisher =     {Kluwer Acad. Publ., Dordrecht},
  series =        {NATO Adv. Sci. Inst. Ser. E: Appl. Sci.},
  title =         {The geometry of bubbles and foams},
  volume =        {354},
  year =          {1999},
}

@article{CoxGranerEtAl,
  author =        {S. J. Cox and F. Graner and F\'Atima Vaz and
                   C. Monnereau-Pittet and N. Pittet},
  journal =       {Philosophical Magazine},
  number =        {11},
  pages =         {1393--1406},
  publisher =     {Taylor \& Francis},
  title =         {Minimal perimeter for N identical bubbles in two
                   dimensions: Calculations and simulations},
  volume =        {83},
  year =          {2003},
  doi =           {10.1080/1478643031000077351},
  url =           {https://doi.org/10.1080/1478643031000077351},
}

@article{CoxMorganGraner,
  author =        {Cox, S. J. and Morgan, F. and Graner, F.},
  journal =       {Proc. R. Soc. Lond. Ser. A Math. Phys. Eng. Sci.},
  number =        {2149},
  pages =         {20120392, 10},
  title =         {Are large perimeter-minimizing two-dimensional
                   clusters of equal-area bubbles hexagonal or
                   circular?},
  volume =        {469},
  year =          {2013},
  doi =           {10.1098/rspa.2012.0392},
  issn =          {1364-5021},
  url =           {https://doi-org.ezlibrary.technion.ac.il/10.1098/
                  rspa.2012.0392},
}

@article{Hales-Honeycomb,
  author =        {Hales, T. C.},
  journal =       {Discrete Comput. Geom.},
  number =        {1},
  pages =         {1--22},
  title =         {The honeycomb conjecture},
  volume =        {25},
  year =          {2001},
  doi =           {10.1007/s004540010071},
  issn =          {0179-5376,1432-0444},
  url =           {https://doi.org/10.1007/s004540010071},
}

@article{HeppesMorgan,
  author =        {Heppes, A. and Morgan, F.},
  journal =       {Phil. Mag.},
  number =        {12},
  pages =         {1333--1345},
  publisher =     {Taylor \& Francis},
  title =         {Planar clusters and perimeter bounds},
  volume =        {85},
  year =          {2005},
  doi =           {10.1080/14786430412331323546},
}

@article{CarocciaMaggi-StableHales,
  author =        {Caroccia, M. and Maggi, F.},
  journal =       {J. Math. Pures Appl. (9)},
  number =        {5},
  pages =         {935--956},
  title =         {A sharp quantitative version of {H}ales'
                   isoperimetric honeycomb theorem},
  volume =        {106},
  year =          {2016},
  doi =           {10.1016/j.matpur.2016.03.017},
  issn =          {0021-7824,1776-3371},
  url =           {https://doi.org/10.1016/j.matpur.2016.03.017},
}

@unpublished{CarocciaDeMasonMaggi-AlmostHoneycomb,
  author =        {Caroccia, M. and DeMason, K. and Maggi, F.},
  note =          {arXiv:2501.05373, to appear in J. Func. Anal.},
  title =         {On the emergence of almost-honeycomb structures in
                   low-energy planar clusters},
  year =          {2025},
}

@unpublished{PeriodicDoubleTilings,
  author =        {Nobili, F. and Novaga, M. and Paolini, E.},
  note =          {arXiv:2502.08396},
  title =         {Periodic double tilings of the plane},
  year =          {2025},
}

@article{OpenProblemsInSoapBubbles96,
  author =        {Sullivan, J. M. and Morgan, F.},
  journal =       {Internat. J. Math.},
  number =        {6},
  pages =         {833--842},
  title =         {Open problems in soap bubble geometry},
  volume =        {7},
  year =          {1996},
  issn =          {0129-167X},
  url =           {https://doi.org/10.1142/S0129167X9600044X},
}

@unpublished{EMilmanXu-Stability,
  author =        {Milman, E. and Xu, B.},
  note =          {arxiv.org:2504.11185},
  title =         {Standard bubbles (and other {M}\"obius-flat
                   partitions) on model spaces are stable},
  year =          {2025},
}

@article{CorneliCorwinEtAl-DoubleBubbleInSandG,
  author =        {Corneli, J. and Corwin, I. and Hurder, S. and
                   Sesum, V. and Xu, Y. and Adams, E. and Davis, D. and
                   Lee, M. and Visocchi, R. and Hoffman, N.},
  journal =       {Houston J. Math.},
  number =        {1},
  pages =         {181--204},
  title =         {Double bubbles in {G}auss space and spheres},
  volume =        {34},
  year =          {2008},
  issn =          {0362-1588},
}

@incollection{Schechtman-ApproxGaussianMultiBubble,
  author =        {Schechtman, G.},
  booktitle =     {Geometric aspects of functional analysis},
  pages =         {373--379},
  publisher =     {Springer, Heidelberg},
  series =        {Lecture Notes in Math.},
  title =         {Approximate {G}aussian isoperimetry for {$k$} sets},
  volume =        {2050},
  year =          {2012},
  url =           {https://doi.org/10.1007/978-3-642-29849-3_23},
}

@article{MontesinosStandardBubbleE!,
  author =        {Montesinos Amilibia, A.},
  journal =       {Asian J. Math.},
  number =        {1},
  pages =         {25--31},
  title =         {Existence and uniqueness of standard bubble clusters
                   of given volumes in {$\Bbb R^N$}},
  volume =        {5},
  year =          {2001},
  doi =           {10.4310/AJM.2001.v5.n1.a3},
  issn =          {1093-6106},
  url =           {https://doi.org/10.4310/AJM.2001.v5.n1.a3},
}

@article{EMilmanNeeman-TripleAndQuadruple,
  author =        {Milman, E. and Neeman, J.},
  journal =       {Acta Math.},
  number =        {1},
  pages =         {71--188},
  title =         {The Structure of Isoperimetric Bubbles on
                   $\mathbb{R}^n$ and $\mathbb{S}^n$},
  volume =        {234},
  year =          {2025},
}

@article{EMilmanNeeman-GaussianMultiBubble,
  author =        {Milman, E. and Neeman, J.},
  journal =       {Ann. of Math. (2)},
  number =        {1},
  pages =         {89--206},
  title =         {The {G}aussian double-bubble and multi-bubble
                   conjectures},
  volume =        {195},
  year =          {2022},
  doi =           {10.4007/annals.2022.195.1.2},
  issn =          {0003-486X},
  url =           {https://doi.org/10.4007/annals.2022.195.1.2},
}

@article{SMALL93,
  author =        {Foisy, J. and Alfaro, M. and Brock, J. and Hodges, N. and
                   Zimba, J.},
  journal =       {Pacific J. Math.},
  number =        {1},
  pages =         {47--59},
  title =         {The standard double soap bubble in {${\bf R}^2$}
                   uniquely minimizes perimeter},
  volume =        {159},
  year =          {1993},
  issn =          {0030-8730},
  url =           {http://projecteuclid.org/euclid.pjm/1102634378},
}

@article{Masters-DoubleBubbleInS2,
  author =        {Masters, J. D.},
  journal =       {Real Anal. Exchange},
  number =        {2},
  pages =         {645--654},
  title =         {The perimeter-minimizing enclosure of two areas in
                   {$S^2$}},
  volume =        {22},
  year =          {1996/97},
  issn =          {0147-1937},
}

@article{CottonFreeman-DoubleBubbleInSandH,
  author =        {Cotton, A. and Freeman, D.},
  journal =       {Int. J. Math. Math. Sci.},
  number =        {11},
  pages =         {641--699},
  title =         {The double bubble problem in spherical space and
                   hyperbolic space},
  volume =        {32},
  year =          {2002},
  issn =          {0161-1712},
  url =           {https://doi.org/10.1155/S0161171202207188},
}

@article{Wichiramala-TripleBubbleInR2,
  author =        {Wichiramala, W.},
  journal =       {J. Reine Angew. Math.},
  pages =         {1--49},
  title =         {Proof of the planar triple bubble conjecture},
  volume =        {567},
  year =          {2004},
  issn =          {0075-4102},
  url =           {https://doi.org/10.1515/crll.2004.011},
}

@article{Lawlor-TripleBubbleInR2AndS2,
  author =        {Lawlor, G. R.},
  journal =       {Anal. Geom. Metr. Spaces},
  number =        {1},
  pages =         {45--61},
  title =         {Perimeter-minimizing triple bubbles in the plane and
                   the 2-sphere},
  volume =        {7},
  year =          {2019},
  doi =           {10.1515/agms-2019-0004},
  url =           {https://doi.org/10.1515/agms-2019-0004},
}

@article{PaoliniTamagnini-PlanarQuadraupleBubbleEqualAreas,
  author =        {Paolini, E. and Tamagnini, A.},
  journal =       {ESAIM Control Optim. Calc. Var.},
  number =        {3},
  pages =         {1303--1331},
  title =         {Minimal clusters of four planar regions with the same
                   area},
  volume =        {24},
  year =          {2018},
  doi =           {10.1051/cocv/2017066},
  issn =          {1292-8119},
  url =           {https://doi.org/10.1051/cocv/2017066},
}

@article{PaoliniTortorelli-PlanarQuadrupleEqualAreas,
  author =        {Paolini, E. and Tortorelli, V. M.},
  journal =       {Calc. Var. Partial Differential Equations},
  number =        {1},
  pages =         {Paper No. 20, 9},
  title =         {The quadruple planar bubble enclosing equal areas is
                   symmetric},
  volume =        {59},
  year =          {2020},
  doi =           {10.1007/s00526-019-1687-9},
  issn =          {0944-2669},
  url =           {https://doi.org/10.1007/s00526-019-1687-9},
}

@article{DoubleBubbleInR3-Announcement,
  author =        {Hutchings, M. and Morgan, F. and Ritore, M. and
                   Ros, A.},
  journal =       {Electron. Res. Announc. Amer. Math. Soc.},
  pages =         {45--49},
  title =         {Proof of the double bubble conjecture},
  volume =        {6},
  year =          {2000},
  doi =           {10.1090/S1079-6762-00-00079-2},
  issn =          {1079-6762},
  url =           {https://doi.org/10.1090/S1079-6762-00-00079-2},
}

@article{DoubleBubbleInR3,
  author =        {Hutchings, M. and Morgan, F. and Ritor\'e, M. and
                   Ros, A.},
  journal =       {Ann. of Math. (2)},
  number =        {2},
  pages =         {459--489},
  title =         {Proof of the double bubble conjecture},
  volume =        {155},
  year =          {2002},
  issn =          {0003-486X},
  url =           {https://doi.org/10.2307/3062123},
}

@article{HHS95,
  author =        {Hass, J. and Hutchings, M. and Schlafly, R.},
  journal =       {Electron. Res. Announc. Amer. Math. Soc.},
  number =        {3},
  pages =         {98--102},
  title =         {The double bubble conjecture},
  volume =        {1},
  year =          {1995},
  issn =          {1079-6762},
  url =           {https://doi.org/10.1090/S1079-6762-95-03001-0},
}

@article{HassSchlafly-EqualDoubleBubbles,
  author =        {Hass, J. and Schlafly, R.},
  journal =       {Ann. of Math. (2)},
  number =        {2},
  pages =         {459--515},
  title =         {Double bubbles minimize},
  volume =        {151},
  year =          {2000},
  doi =           {10.2307/121042},
  issn =          {0003-486X},
  url =           {https://doi.org/10.2307/121042},
}

@article{Hutchings-StructureOfDoubleBubbles,
  author =        {Hutchings, M.},
  journal =       {J. Geom. Anal.},
  number =        {2},
  pages =         {285--304},
  title =         {The structure of area-minimizing double bubbles},
  volume =        {7},
  year =          {1997},
  issn =          {1050-6926},
  url =           {https://doi.org/10.1007/BF02921724},
}

@article{SMALL03,
  author =        {Reichardt, B. W. and Heilmann, C. and Lai, Y. Y. and
                   Spielman, A.},
  journal =       {Pacific J. Math.},
  number =        {2},
  pages =         {347--366},
  title =         {Proof of the double bubble conjecture in {${\bf
                   R}^4$} and certain higher dimensional cases},
  volume =        {208},
  year =          {2003},
  issn =          {0030-8730},
  url =           {https://doi.org/10.2140/pjm.2003.208.347},
}

@article{Reichardt-DoubleBubbleInRn,
  author =        {Reichardt, B. W.},
  journal =       {J. Geom. Anal.},
  number =        {1},
  pages =         {172--191},
  title =         {Proof of the double bubble conjecture in {$\bold
                   R^n$}},
  volume =        {18},
  year =          {2008},
  issn =          {1050-6926},
  url =           {https://doi.org/10.1007/s12220-007-9002-y},
}

@article{Lawlor-DoubleBubbleInRn,
  author =        {Lawlor, G. R.},
  journal =       {J. Geom. Anal.},
  number =        {1},
  pages =         {190--204},
  title =         {Double bubbles for immiscible fluids in
                   {$\Bbb{R}^n$}},
  volume =        {24},
  year =          {2014},
  doi =           {10.1007/s12220-012-9333-1},
  issn =          {1050-6926},
  url =           {https://doi.org/10.1007/s12220-012-9333-1},
}

@article{CorneliHoffmanEtAl-DoubleBubbleIn3D,
  author =        {Corneli, J. and Hoffman, N. and Holt, P. and Lee, G. and
                   Leger, N. and Moseley, S. and Schoenfeld, E.},
  journal =       {J. Geom. Anal.},
  number =        {2},
  pages =         {189--212},
  title =         {Double bubbles in {${\bf S}^3$} and {${\bf H}^3$}},
  volume =        {17},
  year =          {2007},
  issn =          {1050-6926},
  url =           {https://doi.org/10.1007/BF02930720},
}

@unpublished{EMilmanNeeman-QuintupleBubble,
  author =        {Milman, E. and Neeman, J.},
  note =          {arxiv.org:2307.08164},
  title =         {Plateau Bubbles and the Quintuple-Bubble Conjecture},
  year =          {2023},
}

@article{ABV-LensClusters,
  author =        {Alama, S. and Bronsard, L. and Vriend, S.},
  journal =       {Trans. Amer. Math. Soc. Ser. B},
  pages =         {516--535},
  title =         {The standard lens cluster in {$\Bbb {R}^2$} uniquely
                   minimizes relative perimeter},
  volume =        {12},
  year =          {2025},
  doi =           {10.1090/btran/176},
  url =           {https://doi-org.ezlibrary.technion.ac.il/10.1090/btran/176},
}

@article{NPT-LocallyIsoperimetricPartitions,
  author =        {Novaga, M. and Paolini, E. and Tortorelli, V. M.},
  journal =       {Trans. Amer. Math. Soc.},
  number =        {4},
  pages =         {2517--2548},
  title =         {Locally isoperimetric partitions},
  volume =        {378},
  year =          {2025},
  doi =           {10.1090/tran/9339},
  issn =          {0002-9947},
  url =           {https://doi-org.ezlibrary.technion.ac.il/10.1090/tran/9339},
}

@unpublished{BronsardNovak-DifferentTensions,
  author =        {Bronsard, L. and Novack, M.},
  note =          {arXiv:2401.08063, to appear in Ann. Inst. H.
                   Poincar\'e C Anal. Non Lin\'eaire},
  title =         {An Infinite Double Bubble Theorem},
  year =          {2024},
}

@unpublished{BronsardNovackEtAl-NonUniquenessOfLens,
  author =        {Bronsard, L. and Neumayer, R. and Novack, M. and
                   Skorobogatova, A.},
  note =          {arXiv:2507.13995},
  title =         {On the non-uniqueness of locally minimizing clusters
                   via singular cones},
  year =          {2025},
}

@UNPUBLISHED{NPT-NonUniquenessOfLensInR8,
    AUTHOR = {Novaga, M. and Paolini, E. and Tortorelli, V. M.},
    title =        {Existence of a non-standard isoperimetric triple partition},
    note =         {arXiv:2507.14112},
  year =         {2025}
}

@article{CLM-StableDoubleBubbleInR2,
  author =        {Cicalese, M. and Leonardi, G. P. and Maggi, F.},
  journal =       {Interfaces Free Bound.},
  number =        {3},
  pages =         {305--350},
  title =         {Sharp stability inequalities for planar double
                   bubbles},
  volume =        {19},
  year =          {2017},
  doi =           {10.4171/IFB/384},
  issn =          {1463-9963,1463-9971},
  url =           {https://doi.org/10.4171/IFB/384},
}

@article{NPST-ClustersViaConcentrationCompactness,
  author =        {Novaga, M. and Paolini, E. and Stepanov, E. and
                   Tortorelli, V. M.},
  journal =       {J. Geom. Anal.},
  number =        {11},
  pages =         {Paper No. 263, 23},
  title =         {Isoperimetric clusters in homogeneous spaces via
                   concentration compactness},
  volume =        {32},
  year =          {2022},
  doi =           {10.1007/s12220-022-01009-8},
  issn =          {1050-6926,1559-002X},
  url =           {https://doi.org/10.1007/s12220-022-01009-8},
}

@article{NPST-ClustersWithInfinitelyManyCells,
  author =        {Novaga, M. and Paolini, E. and Stepanov, E. and
                   Tortorelli, V. M.},
  journal =       {Netw. Heterog. Media},
  number =        {3},
  pages =         {1226--1235},
  title =         {Isoperimetric planar clusters with infinitely many
                   regions},
  volume =        {18},
  year =          {2023},
  doi =           {10.3934/nhm.2023053},
  issn =          {1556-1801,1556-181X},
  url =           {https://doi.org/10.3934/nhm.2023053},
}

@article{DeRosaTione-ConvexStationaryBubbles,
  author =        {De Rosa, A. and Tione, R.},
  journal =       {Trans. Amer. Math. Soc.},
  number =        {5},
  pages =         {3393--3444},
  title =         {The double and triple bubble problem for stationary
                   varifolds: the convex case},
  volume =        {378},
  year =          {2025},
  doi =           {10.1090/tran/9400},
  issn =          {0002-9947},
  url =           {https://doi-org.ezlibrary.technion.ac.il/10.1090/tran/9400},
}

@article{BCT-StableLensPartitionInR2,
  author =        {Bonacini, M. and Cristoferi, R. and Topaloglu, I.},
  journal =       {Proc. Roy. Soc. Edinburgh Sect. A},
  pages =         {1--34},
  title =         {A stability inequality for planar lens partition},
  year =          {2025},
}

@article{AlmgrenMemoirs,
  author =        {Almgren, Jr., F. J.},
  journal =       {Mem. Amer. Math. Soc.},
  number =        {165},
  title =         {Existence and regularity almost everywhere of
                   solutions to elliptic variational problems with
                   constraints},
  volume =        {4},
  year =          {1976},
}

@article{Taylor-SoapBubbleRegularityInR3,
  author =        {Taylor, J. E.},
  journal =       {Ann. of Math. (2)},
  number =        {3},
  pages =         {489--539},
  title =         {The structure of singularities in soap-bubble-like
                   and soap-film-like minimal surfaces},
  volume =        {103},
  year =          {1976},
  issn =          {0003-486X},
  url =           {https://doi.org/10.2307/1970949},
}

@article{CES-RegularityOfMinimalSurfacesNearCones,
  author =        {Colombo, M. and Edelen, N. and Spolaor, L.},
  journal =       {J. Differential Geom.},
  number =        {3},
  pages =         {411--503},
  title =         {The singular set of minimal surfaces near polyhedral
                   cones},
  volume =        {120},
  year =          {2022},
  doi =           {10.4310/jdg/1649953512},
  issn =          {0022-040X},
  url =           {https://doi.org/10.4310/jdg/1649953512},
}

@incollection{White-AusyAnnouncementOfClusterRegularity,
  author =        {White, B.},
  booktitle =     {Miniconference on geometry and partial differential
                   equations ({C}anberra, 1985)},
  pages =         {244--249},
  publisher =     {Austral. Nat. Univ., Canberra},
  series =        {Proc. Centre Math. Anal. Austral. Nat. Univ.},
  title =         {Regularity of the singular sets in immiscible fluid
                   interfaces and solutions to other {P}lateau-type
                   problems},
  volume =        {10},
  year =          {1986},
}

@article{Simon-Codimension2Regularity,
  author =        {Simon, L.},
  journal =       {J. Differential Geom.},
  number =        {3},
  pages =         {585--652},
  title =         {Cylindrical tangent cones and the singular set of
                   minimal submanifolds},
  volume =        {38},
  year =          {1993},
  issn =          {0022-040X},
  url =           {http://projecteuclid.org/euclid.jdg/1214454484},
}

@article{NaberValtorta-MinimizingHarmonicMaps,
  author =        {Naber, A. and Valtorta, D.},
  journal =       {Ann. of Math. (2)},
  number =        {1},
  pages =         {131--227},
  title =         {Rectifiable-{R}eifenberg and the regularity of
                   stationary and minimizing harmonic maps},
  volume =        {185},
  year =          {2017},
  doi =           {10.4007/annals.2017.185.1.3},
  issn =          {0003-486X},
  url =           {https://doi.org/10.4007/annals.2017.185.1.3},
}

@article{KNS,
  author =        {Kinderlehrer, D. and Nirenberg, L. and Spruck, J.},
  journal =       {J. Analyse Math.},
  pages =         {86--119 (1979)},
  title =         {Regularity in elliptic free boundary problems},
  volume =        {34},
  year =          {1978},
  doi =           {10.1007/BF02790009},
  issn =          {0021-7670},
  url =           {https://doi.org/10.1007/BF02790009},
}

@article{ADN2,
  author =        {Agmon, S. and Douglis, A. and Nirenberg, L.},
  journal =       {Comm. Pure Appl. Math.},
  pages =         {35--92},
  title =         {Estimates near the boundary for solutions of elliptic
                   partial differential equations satisfying general
                   boundary conditions. {II}},
  volume =        {17},
  year =          {1964},
  doi =           {10.1002/cpa.3160170104},
  issn =          {0010-3640},
  url =           {https://doi.org/10.1002/cpa.3160170104},
}

@article{Leonardi-Infiltration,
  author =        {Leonardi, G. P.},
  journal =       {Proc. Roy. Soc. Edinburgh Sect. A},
  number =        {2},
  pages =         {425--436},
  title =         {Infiltrations in immiscible fluids systems},
  volume =        {131},
  year =          {2001},
  doi =           {10.1017/S0308210500000937},
  issn =          {0308-2105},
  url =           {https://doi.org/10.1017/S0308210500000937},
}

@article{HuangMaZhu-ReillyFormula,
  author =        {Huang, G. and Ma, B. and Zhu, M.},
  journal =       {Differential Geom. Appl.},
  pages =         {Paper No. 102136, 20},
  title =         {A {R}eilly type integral formula and its
                   applications},
  volume =        {94},
  year =          {2024},
  doi =           {10.1016/j.difgeo.2024.102136},
  issn =          {0926-2245,1872-6984},
  url =           {https://doi.org/10.1016/j.difgeo.2024.102136},
}

@article{HuangZhu-ConjugatedBrascampLieb,
  author =        {Huang, G. and Zhu, M.},
  journal =       {Proc. Amer. Math. Soc.},
  number =        {11},
  pages =         {4961--4970},
  title =         {Some integral inequalities on weighted {R}iemannian
                   manifolds with boundary},
  volume =        {151},
  year =          {2023},
  doi =           {10.1090/proc/16479},
  issn =          {0002-9939,1088-6826},
  url =           {https://doi.org/10.1090/proc/16479},
}

@incollection{Morgan-StrictCalibrations,
  author =        {Morgan, F.},
  journal =       {Mat. Contemp.},
  note =          {IX School of Differential Geometry (Vit\'oria, 1994)},
  pages =         {139--152},
  title =         {Strict calibrations, constant mean curvature and
                   triple junctions},
  volume =        {9},
  year =          {1995},
  issn =          {0103-9059,2317-6636},
}

@article{LawlorMorgan-TripleJunctions,
  author =        {Lawlor, G. and Morgan, F.},
  journal =       {J. Differential Geom.},
  number =        {3},
  pages =         {514--528},
  title =         {Curvy slicing proves that triple junctions locally
                   minimize area},
  volume =        {44},
  year =          {1996},
  issn =          {0022-040X,1945-743X},
  url =           {http://projecteuclid.org/euclid.jdg/1214459219},
}

@article{McGonagleRoss:15,
  author =        {McGonagle, M. and Ross, J.},
  journal =       {Geom. Dedicata},
  number =        {1},
  pages =         {277--296},
  title =         {The hyperplane is the only stable, smooth solution to
                   the isoperimetric problem in {G}aussian space},
  volume =        {178},
  year =          {2015},
  issn =          {0046-5755},
  url =           {https://doi.org/10.1007/s10711-015-0057-9},
}

@article{Alexandrov-MethodOfMovingPlanes,
  author =        {Aleksandrov, A. D.},
  journal =       {Amer. Math. Soc. Transl. (2)},
  pages =         {412--416},
  title =         {Uniqueness theorems for surfaces in the large. {V}},
  volume =        {21},
  year =          {1962},
  issn =          {0065-9290},
}

@article{BarbosaDoCarmo-StabilityInRn,
  author =        {Barbosa, J. L. and do Carmo, M.},
  journal =       {Math. Z.},
  number =        {3},
  pages =         {339--353},
  title =         {Stability of hypersurfaces with constant mean
                   curvature},
  volume =        {185},
  year =          {1984},
  doi =           {10.1007/BF01215045},
  issn =          {0025-5874},
  url =           {https://doi.org/10.1007/BF01215045},
}

@article{BarbosaDoCarmoEschenburg-StabilityInManifolds,
  author =        {Barbosa, J. L. and do Carmo, M. and Eschenburg, J.},
  journal =       {Math. Z.},
  number =        {1},
  pages =         {123--138},
  title =         {Stability of hypersurfaces of constant mean curvature
                   in {R}iemannian manifolds},
  volume =        {197},
  year =          {1988},
  doi =           {10.1007/BF01161634},
  issn =          {0025-5874},
  url =           {https://doi.org/10.1007/BF01161634},
}

\end{document}